\documentclass[11pt,onecolumn]{IEEEtran}
\usepackage{cite}
\usepackage{amsmath,amssymb,amsfonts}
\usepackage{algorithmic}
\usepackage{graphicx}
\usepackage{dsfont}
\usepackage{pstricks}
\usepackage{afterpage}
\usepackage{comment}
\usepackage{cite}
\allowdisplaybreaks[1]
\newcommand{\diag}{\mbox{diag}}
\newcommand{\eqdef}{\buildrel \triangle \over =}

\usepackage{textcomp}

\begin{document}
\title{Matrix Pencil Based On-Line Computation of Controller Parameters in Dynamic High-Gain Scaling Controllers for Strict-Feedback Systems}
\author{Prashanth Krishnamurthy and Farshad Khorrami
  \thanks{The authors are with the Control/Robotics Research Laboratory (CRRL), Dept. of ECE, NYU Tandon School of Engineering, Brooklyn, NY 11201, USA (e-mails: prashanth.krishnamurthy@nyu.edu, khorrami@nyu.edu).}
  }

\maketitle

\begin{abstract}
  We propose a new matrix pencil based approach for design of state-feedback and output-feedback stabilizing controllers for a general class of uncertain nonlinear strict-feedback-like systems. While the dynamic controller structure is based on the dual dynamic high-gain scaling based approach, we cast the design procedure within a general matrix pencil structure unlike prior results that utilized conservative algebraic bounds of terms arising in Lyapunov inequalities. The design freedoms in the dynamic controller are extracted in terms of generalized eigenvalues associated with matrix pencils formulated to capture the detailed structures (locations of uncertain terms in the system dynamics and their state dependences) of bounds in the Lyapunov analysis. The proposed approach enables efficient computation of non-conservative bounds with reduced algebraic complexity and enhances feasibility of application of the dual dynamic high-gain scaling based control designs. The proposed approach is developed under both the state-feedback and output-feedback cases and the efficacy of the approach is demonstrated through simulation studies on a numerical example.
\end{abstract}

\bstctlcite{IEEEexample:BSTcontrol}

\section{Introduction}

Robust and adaptive control designs for various classes of uncertain nonlinear systems have been addressed in the literature both in state-feedback and output-feedback contexts \cite{KKK95,Isi99,Kha01}. Various specific classes of systems such as feedforward (e.g., \cite{KA01,Tee92,MP96,SJK97b}) and strict-feedback (e.g., \cite{KKK95,KK03,KK04f,YMVX20}) triangular forms have been considered under various sets of assumptions (e.g., including various types of parametric and functional uncertainties, appended dynamics, state and input time delays, input unmodeled dynamics, etc).
One particular control design approach that has been investigated heavily in the literature is based on high-gain scaling using either static or dynamic scaling terms. High gain approaches have been applied for both controller and observer designs.
Classical results in this direction include the adaptive high-gain
controller given in its basic form by $u=-ry,\dot r=y^2$, which is applicable to
minimum-phase systems with relative-degree one  \cite{KS87,SL90,Ilc96}, and static
high-gain scaling based observers \cite{TP94,Kha96,AK99b,Kha08} which introduce observer
gains $r,\ldots,r^n$ with a constant $r$ to obtain semiglobal results.
Combinations of low-gain and high-gain components in control design have also been considered \cite{Lin09b,Lin22}.
State-dependent scaling techniques for control of nonlinear systems are also addressed in \cite{Ito06}.
A combination of a high-gain observer (with the dynamics of the high gain parameter $r$ being of the form of a scalar differential Riccati equation) and a backstepping controller was developed in \cite{Pra03,KKC02}.
A dual observer/controller dynamic high-gain scaling technique was introduced in
\cite{KK02d,KK04f} that combined dynamic scaling based observer and controller structures to address uncertain strict-feedback-like systems (\cite{KKK95,KK03,KK04f,YMVX20}).

The dual dynamic scaling-based control design approach was shown to be a flexible design technique capable of handling 
uncertain terms dependent on all states and uncertain Input-to-State Stable (ISS)
appended dynamics with nonlinear gains from all the system states and the input.
The dynamic high-gain scaling technique provides a unified control design methodology applicable to state-feedback and output-feedback control of strict-feedback (\cite{KK04f,KAP06,KK07f}) and feedforward (\cite{KK04e,KK07h}) systems as well as state-feedback control of nontriangular polynomially-bounded systems (\cite{KK07a}).
The dynamic high-gain scaling based control methodology is also applicable to disturbance attenuation (\cite{KK07h}) as well as control of systems with state and input time delays (\cite{KK10b}) and systems with input unmodeled dynamics \cite{KK13a}.

However, while the dynamic high-gain scaling based methodology has been shown to provide a robust and flexible design approach for a range of system classes, the computation of the design freedoms in the methodology given a specific dynamic system can be challenging due to the algebraic complexity in computing various upper bounds on terms that appear in the Lyapunov analysis as part of the design. These upper bounds primarily capture effects of uncertain terms in the system dynamics in the context of Lyapunov inequalities arising during control design and stability analysis and are utilized in the computation of various design freedoms (such as the scaling parameter dynamics). Due to the algebraic complexity, it is challenging to compute tight bounds for a given specific system and it is often necessary to utilize conservative upper bounds instead. Utilization of conservative upper bounds results in the effective control gains being much larger than required thus hampering practical performance of the closed-loop system. Hence, while the scaled state vector is handled in the Lyapunov analysis in a matrix structure (via coupled Lyapunov inequalities) that provides an elegant formulation, the computation of the upper bounds of uncertain terms often drops down, for algebraic tractability (especially for higher-order systems), to conservative effectively scalar bounds. These conservative scalar bounds do not capture the specific structure of where and how the uncertain terms in the system dynamics appear (i.e., they address the worst-case uncertainties rather than leveraging information of specific state dependence of the uncertain terms in the system dynamics).     

To address the challenges outlined above in application of the dynamic high-gain scaling based control design methodology, we develop a new framework based on a matrix pencil formulation for computation of the constituent design freedoms in the dynamic observer-controller structure. The matrix pencil formulation for computation of design freedoms was introduced in our earlier conference paper \cite{KK20_matrixpencil_strictfeedback} for strict-feedback-like systems and \cite{KK20_matrixpencil_feedforward,KK21_matrixpencil_feedforward} for feedforward-like systems.
The matrix pencil based approach is motivated by the observation that the principal design freedoms are essentially scalars (constant parameters or functions) and the required properties of these design freedoms are naturally formulated in terms of desired forms of Lyapunov inequalities arising in the closed-loop system analysis. Hence, by writing the relevant Lyapunov inequalities in matrix structures, the required properties of the design freedoms can be captured in terms of desired properties of matrix pencils. In particular, it will be seen that the required properties for obtaining non-conservative values of the design freedoms can be formulated in terms of obtaining the smallest possible ``large enough'' or largest possible ``small enough'' values to satisfy certain matrix inequalities. From a matrix pencil viewpoint, the choice of such smallest possible ``large enough'' or largest possible ``small enough'' values is seen to be related to the generalized eigenvalues of the matrix pencils. Matrix pencils which can be viewed most simply in terms of weighted combinations of matrices have been studied over a long history \cite{matrix_pencils} and have been applied in a variety of contexts including system/parameter identification \cite{Hua1988_matrix_pencil,Almunif2020_matrix_pencil}, stability analysis of nonlinear circuits \cite{Riaza2004_matrix_pencil}, and structural methods for design of robust controllers for linear systems \cite{matrix_lecture_notes}.

The proposed methodology in this paper is based on our earlier conference paper \cite{KK20_matrixpencil_strictfeedback} and builds further through a more detailed development of the methodology including both the state-feedback and output-feedback cases and simulation-based studies of the methodology.
By capturing the specific state dependence structures of the known and unknown terms in the Lyapunov inequalities arising during the control design, the proposed approach enables efficient non-conservative computation of the design freedoms with reduced algebraic complexity. The overall computation of the design freedoms is structured as a sequence of matrix pencil based subproblems that capture the specific properties that each design freedom needs to satisfy in the context of the detailed structure of Lyapunov inequalities that appear in the control design.
By characterizing the required properties of design freedoms in terms of matrix pencil characterizations that can be evaluated at run-time, the design freedoms can be computed on-line instead of requiring pre-computation of conservative worst-case values.
The proposed approach enables efficient computation of less conservative values for the design freedoms and removes the need for algebraically computing the required bounds; both these benefits facilitate application of the dynamic high-gain scaling based control design techniques.   

The structure of this paper is as follows.
The class of nonlinear uncertain systems considered, assumptions imposed on the system structure, and the control design objective are defined in Section~\ref{sec:formulation}. The dynamic high-gain scaling methodology applied to this class of systems is summarized in Section~\ref{sec:scaling_control} under the state-feedback (Section~\ref{sec:scaling_control_sf}) and output-feedback (Section~\ref{sec:scaling_control_of}) cases. The proposed matrix pencil based approach for computation of the design freedoms in the dynamic high-gain scaling based controller is presented in Section~\ref{sec:design_sf} under the state-feedback case. As part of Section~\ref{sec:design_sf}, the application of the proposed approach to an illustrative example is presented in Section~\ref{sec:example_sf}. The application of the proposed matrix pencil based approach under the output-feedback case is addressed in Section~\ref{sec:design_of} including an illustrative numerical example in Section~\ref{sec:example_of}.
Concluding remarks are provided in Section~\ref{sec:conclusion}.

\section{Problem Formulation}
\label{sec:formulation}
We consider nonlinear uncertain strict-feedback-like systems of the form:  
\begin{align}
\dot x_i&=\phi_i(\overline x_i)+\phi_{(i,i+1)}(x_1)x_{i+1} \ \ , i=1,\ldots,n-1 \nonumber\\
  \dot x_n&=\phi_n(x) +\mu_0(x_1)u.
            \label{system}
\end{align}
We consider both the state-feedback case ($y=x$) and the output-feedback case
\begin{align}
  y &= x_1.
      \label{system_output}
\end{align}
In \eqref{system} and \eqref{system_output}, $x=[x_1,\ldots,x_n]^T\in{\mathcal R}^n$ is the state of the system\footnote{${\mathcal R}$, ${\mathcal R}^+$, and ${\mathcal R}^k$ denote the set of real numbers, the set of non-negative real numbers, and the set of real $k$-dimensional column vectors, respectively.}, $u\in{\mathcal R}$ is the input, $\overline x_i$ denotes $[x_1,\ldots,x_i]^T$, and $y\in{\mathcal R}^{n_y}$ is the measured output ($n_y=n$ in the state-feedback case and $n_y=1$ in the output-feedback case). 
In \eqref{system}, $\phi_{(i,i+1)}, i=1,\ldots,n-1$, and $\mu_0$ are known scalar real-valued
continuous functions of the measured output $y$.  
$\phi_i, i=1,\ldots, n$, are
scalar real-valued uncertain
functions of their arguments\footnote{While the arguments for each $\phi_i$ can include the entire state $x$ as long as the bounds in Assumption~A2 below is satisfied, a triangular structure for state dependence of the $\phi_i$ terms is shown in \eqref{system} to highlight the triangular state dependence structure of the bounds to be introduced in Assumption~A2.}.
In the output-feedback case, only the output $y=x_1$ is assumed to be measured and the state variables $x_2,\ldots,x_n$ are unmeasured.
In the state-feedback case, all the state variables $x_1,\ldots,x_n$ are assumed to be measured.
The control design problem being addressed is to design an asymptotically stabilizing dynamic state-feedback/output-feedback control law for $u$ using the measurement of $y$.
The assumptions imposed on the system \eqref{system} are listed below.

\noindent{\bf Assumption A1} ({\em lower boundedness away from zero of ``upper diagonal'' terms $\phi_{(i,i+1)}$ and $\mu_0$}):
A constant $\sigma > 0$ exists such that\footnote{  The notation $|a|$ is used to denote the Euclidean norm of a vector $a$.  If $a$ is a scalar, $|a|$ denotes its absolute value. }
 $|\phi_{(i,i+1)}(x_1)|\geq \sigma ,1\leq i\leq n-1$, and $|\mu_0(x_1)|\geq 0$ for all $x_1\in{\mathcal R}$. Since $\phi_{(i,i+1)}$ and $\mu_0$ are continuous functions, this assumption can, without loss of generality,
 be stated as
 $\phi_{(i,i+1)}(x_1)\geq \sigma>0 ,1\leq i\leq n-1$, and $\mu_0(x_1)\geq \sigma > 0$.

 \noindent{\bf Assumption A2} ({\em Bounds on uncertain functions $\phi_i$}):
The functions $\phi_i, i=1,\ldots,n$, can be bounded as
$|\phi_i(\overline x_i)|\leq
\sum_{j=1}^i \phi_{(i,j)}(x_1)|x_j|$
for all $x\in{\cal R}^n$ where 
$\phi_{(i,j)}(x_1), i=1,\ldots,n, j=1,\ldots,i$,
are known continuous
non-negative functions.

\noindent{\bf Assumption A3} ({\em Bi-directional cascading dominance of ``upper diagonal'' terms $\phi_{(i,i+1)}, i=2,\ldots,n-1$,}):
Positive constants $\overline\rho_i,i=3,\ldots,n-1$, and $\underline\rho_i,i=3,\ldots,n-1$ exist such that
$\forall x_1\in {\cal R}$
\begin{align}
  \phi_{(i,i+1)}(x_1)\geq \overline\rho_i \phi_{(i-1,i)}(x_1)\,,\,\,\, i=3,\ldots,n-1
  \label{eq:cascading_dominance_controller}
  \\
  \phi_{(i,i+1)}(x_1)\leq \underline\rho_i \phi_{(i-1,i)}(x_1)\,,\,\,\, i=3,\ldots,n-1.
  \label{eq:cascading_dominance_observer}
\end{align}

\noindent{\bf Remark 1:}
The structure of the assumptions above are analogous to the assumptions introduced for the dual dynamic high gain based state-feedback and output-feedback control designs in \cite{KK04f} except that Assumption~A2 involves a more detailed structure of the state dependence of each $\phi_i$. While 
Assumption~A2 considers  $|\phi_i(\overline x_i)|\leq \sum_{j=1}^i \phi_{(i,j)}(x_1)|x_j|$ where functions $\phi_{(i,j)}$ model the detailed structure of where and how the uncertainties appear, the
output-feedback control design in \cite{KK04f} considered a bound of the form $|\phi_i(\overline x_i)|\leq \Gamma(x_1)\sum_{j=1}^i |x_j|$, which can be viewed as a ``worst-case'' bound with $\Gamma$ being a single known function instead of separate $\phi_{(i,j)}(x_1)$. The 
Assumption A1 ensures observability, controllability, and uniform relative degree of the system (with $u$ being the control input and $x_1$ the output).
The bounds in Assumption~A2 on uncertain terms $\phi_i$ essentially requires the uncertain terms to be bounded linearly in unmeasured state variables with a triangular state dependence structure in the bounds.
While the bounds can be extended to include parametric uncertainties (without requiring known upper bounds on unknown parameters) via an unknown parameter $\theta$ multiplying the known upper bounds, such unknown parameters are not included in Assumption~A2 for brevity in this paper and to focus on the main aspects of the proposed matrix pencil based approach.
The bounds on relative ``sizes'' (in a nonlinear function sense) of the upper diagonal terms $\phi_{(i,i+1)}$ in
Assumption A3 essentially require these nonlinear functions to be comparable (up to constant factors, e.g., $1+x_1^2$ vs. $1+x_1+x_1^2$) and is vital in achieving solvability of pairs of coupled Lyapunov inequalities in Sections~\ref{sec:coupled_lyap_sf} and \ref{sec:coupled_lyap_of}.
The first set of inequalities in Assumption~A3 given by \eqref{eq:cascading_dominance_controller} is the controller-context cascading dominance condition which essentially requires upper diagonal terms closer to the input to be larger (in a nonlinear function sense).
The second set of inequalities in Assumption~A3 given by \eqref{eq:cascading_dominance_observer} is the observer-context cascading dominance condition which requires upper diagonal terms closer to output $y=x_1$ to be larger.
While both the controller-context and observer-context cascading dominance conditions are written as part of a single Assumption for convenience, it is to be noted that the control design in the state-feedback case only requires controller-context cascading dominance. The controller and observer designs in the output-feedback case require controller-context and observer-context cascading dominance conditions, respectively.
The functions $\phi_{(i,i+1)}$ are referred to here as ``upper diagonal'' terms since if the dynamics \eqref{system} were to be written in the form $\dot x = A(x_1)x + B(x_1)u + \phi(x)$ with $\phi=[\phi_1,\ldots,\phi_n]^T$, the functions $\phi_{(i,i+1)}$ would appear on the upper diagonal of the matrix $A(x_1)$. The upper diagonal terms $\phi_{(i,i+1)}$ and the nonlinear functions $\phi_{(i,j)}$ appearing in the bounds in Assumption~A2 can be allowed to be functions of $y=x$ in the state-feedback case as in \cite{KK04f}; however, these functions are considered in this paper to be dependent only on $x_1$ for algebraic simplicity and to focus on the main aspects of the proposed methodology.

\section{Dual Dynamic Scaling-Based Control Design}
\label{sec:scaling_control}
The dynamic scaling-based control design methodology from \cite{KK04f} introduces a dynamic high-gain parameter and uses successive powers of this parameter to scale the state variables of the system following which the dynamics of the high-gain parameter are designed via a Lyapunov analysis. In the output-feedback case, the high-gain parameter is also used in the definition of the observer gains and scaling of the observer error state variables. The basic structure of the dynamic scaling-based control design from \cite{KK04f} is summarized below in the state-feedback (Section~\ref{sec:scaling_control_sf}) and output-feedback (Section~\ref{sec:scaling_control_of}) cases. 
As seen below, both the observer and controller involve powers of a scaling parameter $r$, which is generated through a dynamics that is designed taking into account the various uncertainties in the system structure.

\subsection{Control Design Under State-Feedback Case}
\label{sec:scaling_control_sf}
The state scaling and structure of the control law are presented in Section~\ref{sec:controller_design_sf}. The design of the controller gains using a pair of coupled Lyapunov inequalities is summarized in Section~\ref{sec:coupled_lyap_sf}. The design freedoms in the controller and the considerations for picking appropriate values for the design freedoms are discussed in Section~\ref{sec:design_freedoms_sf}.

\subsubsection{Controller Design}
\label{sec:controller_design_sf}
Scaled versions of the state variables of the system \eqref{system} are defined as
$\eta=[\eta_2,\ldots,\eta_n]^T$ where $\eta_2,\ldots,\eta_n$ are given by
\begin{align}
  \eta_2 &= \frac{x_2 + \zeta(x_1)}{r}
           \ \ ; \ \ 
  \eta_i = \frac{x_i}{r^{i-1}}, i=3,\ldots,n.
           \label{etaidefn_sf}
\end{align}
In \eqref{etaidefn_sf}, $r$ is a dynamic high-gain scaling parameter, the dynamics of which will be designed later in this section. The dynamics of the high-gain scaling parameter $r$ to be designed will
be such that $r(t)\geq 1$ for all time $t\geq 0$.
Also, in \eqref{etaidefn_sf}, $\zeta(x_1)$ is a function defined to be of the form
\begin{align}
  \zeta(x_1) &= x_1\zeta_1(x_1)
               \label{zetadefn}
\end{align}
where $\zeta_1$ is a function that will be designed further below.
The dynamics of the scaled state vector $\eta$ defined in \eqref{etaidefn_sf} are given by\footnote{For notational convenience, we drop arguments of functions when no confusion will result.} 
\begin{align}
  \dot\eta_2&= r\phi_{(2,3)}\eta_3 + \frac{1}{r}\phi_2
              +\frac{[\zeta_1'(x_1)x_1+\zeta_1(x_1)]}{r}[(r\eta_2-\zeta)\phi_{(1,2)}+\phi_1]
               -\frac{\dot r}{r}\eta_2
              \nonumber\\
  \dot \eta_i&=r\phi_{(i,i+1)}\eta_{i+1} +\frac{1}{r^{i-1}}\phi_i -\frac{\dot r}{r}(i-1)\eta_i \,\, ,\,\, i=3,\ldots,n-1
              \nonumber\\
  \dot \eta_n&= -\frac{\dot r}{r}(n-1)\eta_n +\frac{1}{r^{n-1}}\phi_n +\frac{1}{r^{n-1}}\mu_{0}u
               \label{etadyn_sf}
\end{align}
where $\zeta_1'(x_1)$ denotes the partial derivative 
of $\zeta_1$ with respect to its argument evaluated at $x_1$.
The control input $u$ is designed as
\begin{align}
u&=-\frac{r^n}{\mu_0(x_1)}K_c \eta
\label{sfcontrol}
\end{align}
with $K_c=[k_2,\ldots,k_n]$ where $k_i,i=2,\ldots,n$, are functions of $x_1$. These functions $k_i$, which can be regarded as controller gains appearing in the definition of the control law \eqref{sfcontrol} will be designed in Section~\ref{sec:coupled_lyap_sf} based on a pair of coupled Lyapunov inequalities involving the upper diagonal terms $\phi_{(i,i+1)}$.
The dynamics of $\eta$ under the control law
\eqref{sfcontrol} are given by
\begin{align}
\dot\eta&=rA_c \eta -\frac{\dot r}{r}D_c\eta+ \Phi +H\eta_2+\Xi
\label{sfetamatdyn_sf}
\end{align}
where 
$A_c$  is the $(n-1)\times (n-1)$ matrix with $(i,j)^{th}$
element
\begin{align}
  A_{c_{(i,i+1)}}(x_1)&=\phi_{(i+1,i+2)}(x_1)\,\, ,\,\, i=1,\ldots, n-2
\nonumber\\ 
A_{c_{(n-1,j)}}(x_1)&=-k_{j+1}(x_1)\,\, ,\,\,  j=1,\ldots,n-1
\label{Acdefn}
\end{align}
with zeros elsewhere, and\footnote{The notation $\mbox{diag}(T_1,\ldots,T_m)$ denotes an $m\times m$ diagonal
  matrix with diagonal elements $T_1,\ldots,T_m$.
  If $T_1,\ldots,T_m$ are matrices, the notation $\mbox{diag}(T_1,\ldots,T_m)$ denotes the block diagonal matrix with the blocks on the principal diagonal being $T_1,\ldots,T_m$.
  Also, $\mbox{lowerdiag}(T_1,\ldots,T_{m-1})$ and $\mbox{upperdiag}(T_1,\ldots,T_{m-1})$ denote the $m\times m$ matrices with the lower diagonal entries (i.e., $(i+1,i)^{th}$ entries, $i=1,\ldots,m-1$) and upper diagonal entries (i.e., $(i,i+1)^{th}$ entries, $i=1,\ldots,m-1$), respectively, being $T_1,\ldots,T_{m-1}$ and zeros elsewhere.
  $I_m$ denotes the $m\times m$
  identity matrix.}
\begin{align}
D_c&=\mbox{diag}(1,2,\ldots,n-1)
\label{Dcdefn}
\\
\Phi &= \left[\frac{\phi_2}{r},\ldots,\frac{\phi_n}{r^{n-1}}\right]^T
      \label{Phidefn}
      \\
  H&=[(\zeta_1'(x_1)x_1+\zeta_1)\phi_{(1,2)},0,\ldots,0]^T
     \label{H_defn}
\\
  \Xi&=[\frac{(\phi_1-\zeta\phi_{(1,2)})[\zeta_1'(x_1)x_1+\zeta_1]}{r},
       0,\ldots,0]^T.
       \label{Xi_defn}
\end{align}
With $P_c$ being a symmetric positive definite matrix designed in Section~\ref{sec:coupled_lyap_sf} based on a pair of coupled Lyapunov inequalities, a Lyapunov function is defined for the closed-loop system as
\begin{align}
  V &=  \frac{1}{2}x_1^2 + r\eta^T P_c\eta.
        \label{Vdefn_sf}
\end{align}
From \eqref{sfetamatdyn_sf} and \eqref{Vdefn_sf}, we have
\begin{align}
 \dot V
  &=
    x_1[\phi_1+(r\eta_2 -\zeta_1x_1)\phi_{(1,2)}]
                 +r^2\eta^T[P_c A_c +A_c^T P_c]\eta
  \nonumber\\&\quad 
  +2r\eta^T P_c(\Phi+H\eta_2+\Xi)
  -\dot r\eta^T[P_c \tilde D_c\!+\!\tilde D_c P_c]\eta
  \label{Vdot1_sf}
\end{align}
where
$\tilde D_c=D_c- \frac{1}{2}I_{n-1}$.

\subsubsection{Coupled Lyapunov Inequalities}
\label{sec:coupled_lyap_sf}
The conditions in Assumption A3 on the relative ``sizes'' (in a nonlinear function sense) of the upper diagonal terms $\phi_{(i,i+1)}$ are the {\em cascading dominance} conditions introduced in \cite{KK04f}. The two sets of inequalities in Assumption~A3 are the controller-context cascading dominance conditions in \eqref{eq:cascading_dominance_controller} and the observer-context cascading dominance conditions in \eqref{eq:cascading_dominance_observer}, which are relevant in the designs of the controller and observer, respectively. For the state-feedback case considered in this section, only the controller-context cascading dominance conditions are required. These cascading dominance conditions were shown in \cite{KK04f,KK06} to be closely related to solvability of pairs of coupled Lyapunov inequalities that appear in the high gain based control design.
Under Assumption A1 and condition \eqref{eq:cascading_dominance_controller} in Assumption A3, it is possible to construct (\cite{KK04f,KK06}) a symmetric positive definite matrix $P_c$ and functions $k_2(x_1),\ldots,k_n(x_1)$ (that appear in the definition of matrix $A_c$) such that the following coupled Lyapunov inequalities are satisfied (for all $x_1\in{\cal R}$) with some positive constants $\nu_{c}$, $\underline\nu_{c}$,
and $\overline\nu_{c}$:
\begin{align}
  P_c A_c+A_c^T P_c \leq -\nu_c \phi_{(2,3)} I
  \ \ ; \ \ 
  \underline \nu_c I \leq P_c \tilde D_c + \tilde D_c P_c\leq \overline\nu_c I.
  \label{coupled_lyap}
\end{align}
From Theorem 2 in \cite{KK06}, the functions $k_2,\ldots,k_n$ can be chosen to be linear constant-coefficient combinations of $\phi_{(2,3)},\ldots,\phi_{(n-1,n)}$.

\subsubsection{Design Freedoms}
\label{sec:design_freedoms_sf}
The design freedoms appearing in the Lyapunov inequality \eqref{Vdot1_sf} are the function $\zeta_1$ and the dynamics of $r$. The strategies for picking these design freedoms and the typical high-gain scaling based computations for these design freedoms are outlined below. The basic strategy in dynamic scaling-based designs is to design the dynamics of $r$ in such a way that the derivative $\dot r$ is ``large'' until the scaling parameter itself becomes ``large.'' For this purpose, the dynamics of $r$ can be designed to be of the form 
\begin{align}
  \dot r &= \max\{-ar(r-1) + r \Omega(x_1),0\} \ \ ; \ r(0) \geq 1
  \label{rdot}
\end{align}
where $a$ is to be picked as a positive constant (or positive function of $x_1$ lower bounded by a positive constant) and $\Omega(x_1)$ is a function to be picked taking into account the various terms appearing in \eqref{Vdot1_sf}. When $r$ becomes ``large enough'' (specifically, when $r\geq R(x_1)$ where $R(x_1)$ is defined as $R(x_1)=1 + \frac{\Omega(x_1)}{a}$), we have $\dot r = 0$.
With this design of the dynamics of $r$, \eqref{Vdot1_sf} yields
\begin{align}
  \dot V
  &=
    r^2\eta^T[P_c A_c +A_c^T P_c]\eta
    +a r^2\eta^T[P_c \tilde D_c+\tilde D_c P_c]\eta
  \nonumber\\&\quad 
  + rx_1\eta_2\phi_{(1,2)}
  +2r\eta^T P_c(\Phi+H\eta_2+\Xi)
  \nonumber\\&\quad 
  -r(\Omega + a)\eta^T[P_c \tilde D_c+\tilde D_c P_c]\eta
  \nonumber\\&\quad 
  -\zeta_1 x_1^2 + x_1\phi_1
  \label{Vdot2_sf}
\end{align}
Using \eqref{coupled_lyap}, \eqref{Vdot2_sf} yields
\begin{align}
  \dot V
  &\leq
    -r^2(\nu_c\phi_{(2,3)}-a\overline\nu_c)|\eta|^2
  \nonumber\\&\quad 
  + rx_1\eta_2\phi_{(1,2)}
  +2r\eta^T P_c(\Phi+H\eta_2+\Xi)
  \nonumber\\&\quad 
  -r(\Omega + a)\eta^T[P_c \tilde D_c+\tilde D_c P_c]\eta
  \nonumber\\&\quad 
  -\zeta_1 x_1^2 + x_1\phi_1
  \label{Vdot3_sf}
\end{align}

The terms in \eqref{Vdot2_sf} and \eqref{Vdot3_sf} can be considered in terms of three parts:
    (a) terms multiplied by $r^2$, i.e., the first line of \eqref{Vdot2_sf} or \eqref{Vdot3_sf}; (b)
    terms multiplied by $r$, i.e., the second and third lines of \eqref{Vdot2_sf} or \eqref{Vdot3_sf}; (c)
    terms not explicitly multiplied by $r$, i.e., the last line of \eqref{Vdot2_sf} or \eqref{Vdot3_sf}.

Based on \eqref{Vdot3_sf}, the primary design considerations are summarized below:
\begin{itemize}
\item The positive quantity $a$ appearing in the dynamics of $r$ in \eqref{rdot} must be picked small enough to retain the negativeness of the terms multiplied by $r^2$ in the  first line of \eqref{Vdot3_sf}. 
\item The functions $\zeta_1$ and $\Omega_1$ must be picked such that the sign-indefinite terms in the second line and the second term in the last line of \eqref{Vdot3_sf} are dominated by the terms in the third line and the first term in the last line.
\end{itemize}
Based on these design considerations, $a$ can, for example, be picked such that
  \begin{align}
    a &\leq \frac{\nu_c\sigma}{2\overline\nu_c}.
        \label{a_choice1_sf}
  \end{align}
  The design of $\zeta_1$ and $\Omega$ entails computing upper bounds on the sign-indefinite terms in the second line of \eqref{Vdot3_sf} and the second term in the last line of \eqref{Vdot3_sf}. For example,
  using the bounds on the uncertain terms $\phi_i$ in Assumption~A2 and the property that $r\geq 1$, an upper bound for $\Phi$ can be computed as\footnote{Given a matrix $M$, the notation $||M||$ denotes its Frobenius norm.}
  \begin{align}
    |\Phi| &\leq \frac{|x_1|}{r}[|\tilde\phi_1| + |\zeta_1(x_1)||\tilde\phi_2|]
                      + ||\tilde A(x_1)|| |\eta|
                      \label{eq:Phibarbound_sf}
  \end{align}
  where
  \begin{align}
  \tilde\phi_1 &= [\phi_{(2,1)},\phi_{(3,1)},\ldots,\phi_{(n,1)}]^T \ \ ; \ \ 
                 \tilde\phi_2 = [\phi_{(2,2)},\phi_{(3,2)},\ldots,\phi_{(n,2)}]^T
                 \label{tilde_phi1_phi2_defn}
  \end{align}
  and
  $\tilde A$ denotes the $(n-1)\times (n-1)$ matrix with $(i,j)^{th}$ element $\phi_{(i+1,j+1)}$ at locations on and below the diagonal and zeros everywhere else. From \eqref{eq:Phibarbound_sf}, an upper bound on
$2r\eta^T P_c\Phi$ can be computed. Similarly, a term such as $2r\eta^T P_cH\eta_2$ can be upper bounded as\footnote{Given a symmetric positive-definite matrix $P$, 
    $\lambda_{max}(P)$ and $\lambda_{min}(P)$ denote its maximum and minimum eigenvalues, respectively.} 
  \begin{align}
    2r\eta^T P_cH\eta_2 &\leq 3r \lambda_{max}(P_c) \phi_{(1,2)}|\zeta_1'x_1+\zeta_1|
      |\eta|^2.
      \label{bnd_ineq1_sf}
  \end{align}
  Also,
  \begin{align}
    x_1\phi_1 &\leq  x_1^2 \phi_{(1,1)}(x_1).
                \label{bnd_ineq2_sf}
  \end{align}
  Based on similarly derived upper bounds for other terms in the second line of \eqref{Vdot3_sf}, the functions $\zeta_1$ and $\Omega$ can then be designed so that the terms in the fourth line and first term of the last line of \eqref{Vdot3_sf} dominate the obtained upper bounds.  

While the typical design procedure in dynamic scaling-based control designs (and indeed in most other nonlinear control designs) is based on the sort of algebraic computation of upper bounds as described above, this is algebraically complex and often results in conservative upper bounds for algebraic tractability (for instance, the bounds in \eqref{eq:Phibarbound_sf} and \eqref{bnd_ineq1_sf} are essentially worst-case upper bounds that do not, for example, take into account the relative sizes of functions $\phi_{(i,j)}$). In contrast, the matrix pencil based approach described in Section~\ref{sec:design_sf} will formulate the choice of the design freedoms discussed above directly in terms of the required properties to be enforced on the Lyapunov inequalities to thereby obtain non-conservative bounds with lower algebraic complexity. 

\subsection{Control Design Under Output-Feedback Case}
\label{sec:scaling_control_of}
In a  structure parallel to the development of the state-feedback controller in Section~\ref{sec:scaling_control_sf}, the dual dynamic scaling-based control design under the output-feedback case is considered below and the observer and controller structures, corresponding coupled Lyapunov inequalities, and design freedoms are discussed in Sections~\ref{sec:obs_controller_design_of}, \ref{sec:coupled_lyap_of}, and \ref{sec:design_freedoms_of}, respectively.

\subsubsection{Observer and Controller Designs}
\label{sec:obs_controller_design_of}
A  dynamic reduced-order observer is applied of the form
\begin{align}
  \dot{\hat x}_i &= \phi_{(i,i+1)}(x_1)[\hat x_{i+1} + r^i f_{i+1}(x_1)]
                   - r^{i-1}g_i(x_1)[\hat x_2 + rf_2(x_1)] 
  \nonumber\\
  &\quad - (i-1)\dot r r^{i-2}f_i(x_1), 2\leq i \leq n-1
  \nonumber\\
  \dot{\hat x}_n &= \mu_{0}(x_1) u 
                   - r^{n-1}g_n(x_1)[\hat x_2 + rf_2(x_1)]
  - (n-1)\dot r r^{n-2}f_n(x_1).
  \label{eq:obs}
\end{align}
where
$\hat x = [\hat x_2,\ldots, \hat x_n]^T\in{\mathcal R}^{n-1}$ is the state of the reduced-order observer.
In \eqref{eq:obs}, $g_i(x_1)$ are observer gain functions that will be designed in Section~\ref{sec:coupled_lyap_of} based on a pair of coupled Lyapunov inequalities involving the upper diagonal terms $\phi_{(i,i+1)}$ analogous to the design of the controller gain functions in Section~\ref{sec:coupled_lyap_sf}.
Also, $f_i(x_1)$ are functions defined as
    \begin{align}
      f_i(x_1) &= \int_0^{x_1}\frac{g_i(\pi)}{\phi_{(1,2)}(\pi)}d\pi, 2\leq i \leq n.
      \label{eq:f_i}
    \end{align}
    As in the state-feedback case in Section~\ref{sec:scaling_control_sf},
 the dynamics of the high-gain scaling parameter $r$ to be designed will
be such that $r(t)\geq 1$ for all time $t\geq 0$.

In the reduced-order observer dynamics defined above, the estimates for the unmeasured states $x_i$ are $\hat x_i + r^{i-1}f_i(x_1)$ and the observer errors are defined as
\begin{align}
  e_i &= \hat x_i + r^{i-1}f_i(x_1) - x_i, 2\leq i \leq n
  \label{eidefn}
\end{align}
and the scaled observer errors (scaled by powers of the high-gain scaling parameter $r$) are defined as
\begin{align}
  \epsilon_i &= \frac{e_i}{r^{i-1}}, i=2\leq i \leq n
  \ \ ; \ 
  \epsilon = [\epsilon_2,\ldots,\epsilon_n]^T.
\end{align}
The dynamics of $\epsilon$ can be written in a matrix structure as
\begin{align}
 \dot\epsilon &= r A_o \epsilon - \frac{\dot r}{r}D_o\epsilon +\overline\Phi 
 \label{sfobs}
\end{align}
where
\begin{itemize}
  \item
$A_o$ is a $(n-1)\times (n-1)$ matrix with $(i,j)^{th}$ entry
\begin{align}
  A_{o_{(i,1)}}(x_1) &= -g_{i+1}(x_1) \ \ , \  i=1,\ldots,n-1
  \nonumber\\
  A_{o_{i,i+1}}(x_1) &= \phi_{(i+1,i+2)}(x_1) \ \ , \ i=1,\ldots,n-2
\end{align}
and with zeros everywhere else
\item $D_o$ is a $(n-1)\times (n-1)$ matrix defined as
 $D_o=\mbox{diag}(1,2,\ldots,n-1)$
  \item $\overline\Phi$ is given by
\end{itemize}
\begin{align}
  \!\!\!\! \overline\Phi &= [\overline\Phi_2,\ldots,\overline\Phi_n]^T
  ; 
  \overline\Phi_i = -\frac{\phi_i(\overline x_i)}{r^{i-1}} \!+\! g_i(x_1)\frac{\phi_1(x_1)}{\phi_{(1,2)}(x_1)}.
                    \label{Phiidefn}
\end{align}
The scaled estimates $\hat x_i + r^{i-1}f_i(x_1)$ for the unmeasured state variables $x_i, i=2,\ldots,n$, are defined as
$\eta=[\eta_2,\ldots,\eta_n]^T$ where $\eta_2,\ldots,\eta_n$ are given by
\begin{align}
  \eta_2 &= \frac{\hat x_2 + rf_2(x_1) + \zeta(x_1)}{r}
           \ \ ; \ \ 
  \eta_i = \frac{\hat x_i + r^{i-1}f_i(x_1)}{r^{i-1}}, i=3,\ldots,n.
           \label{etaidefn_of}
\end{align}
The definitions of $\eta_2,\ldots,\eta_n$ are analogous to the corresponding definitions in \eqref{etaidefn_sf} in the state-feedback case with the unmeasured state variables $x_2,\ldots,x_n$ replaced with their corresponding observer-based estimates $\hat x_2+rf_2,\ldots,\hat x_n+r^{n-1}f_n$.
Also, analogous to \eqref{etaidefn_sf}, $\zeta(x_1)$ in \eqref{etaidefn_of} is a function defined to be of the form shown in \eqref{zetadefn}
with $\zeta_1$ being a function to be designed.
The dynamics of the scaled observer estimates $\eta_i$ defined in \eqref{etaidefn_of} are given by
\begin{align}
\dot\eta_2&=
            r\phi_{(2,3)}\eta_3 - rg_2\epsilon_2
          +\frac{[\zeta_1'(x_1)x_1+\zeta_1(x_1)]}{r}[(r\eta_2-\zeta-r\epsilon_2)\phi_{(1,2)}+\phi_1]
            \nonumber\\
           &\quad +g_2\frac{\phi_1}{\phi_{(1,2)}}  -\frac{\dot r}{r}\eta_2
\nonumber\\
\dot \eta_i&=r\phi_{(i,i+1)}\eta_{i+1}-rg_i\epsilon_2+g_i\frac{\phi_1}{\phi_{(1,2)}}
 -\frac{\dot r}{r}(i-1)\eta_i \,\, ,\,\, i=3,\ldots,n-1
             \nonumber\\
\dot \eta_n&= -rg_n\epsilon_2 + g_n\frac{\phi_1}{\phi_{(1,2)}}-\frac{\dot r}{r}(n-1)\eta_n +\frac{1}{r^{n-1}}\mu_{0}u.
\end{align}

The control input $u$ is designed in terms of the vector $\eta$ of scaled observer estimates defined above as \eqref{sfcontrol} analogous to the state-feedback case
with $K_c=[k_2,\ldots,k_n]$ being the controller gain vector where $k_i,i=2,\ldots,n$, are functions of $x_1$ designed as in Section~\ref{sec:coupled_lyap_sf} based on a pair of coupled Lyapunov inequalities.
The dynamics of $\eta$ under the control law
\eqref{sfcontrol} are given by
\begin{align}
\dot\eta&=rA_c \eta -\frac{\dot r}{r}D_c\eta+\Phi -rG\epsilon_2 +H[\eta_2-\epsilon_2]+\Xi
\label{sfetamatdyn_of}
\end{align}
where $A_c$, $D_c$, $H$, and $\Xi$ are the same as in the corresponding dynamics \eqref{sfetamatdyn_sf} in the state-feedback case (equations \eqref{Acdefn}, \eqref{Dcdefn}, \eqref{H_defn}, and \eqref{Xi_defn} respectively), and $G$ and $\Phi$ are given by
\begin{align}
      G &= [g_2,\ldots,g_n]^T
\ \ ; \ \ 
\Phi=\frac{\phi_1}{\phi_{(1,2)}}G.
      \label{G_Phidefn}
\end{align}
With $P_c$ and $P_o$ being symmetric positive definite matrices designed as in Sections~\ref{sec:coupled_lyap_sf} and \ref{sec:coupled_lyap_of}, respectively, the observer and controller Lyapunov functions $V_o$ and $V_c$ and the overall Lyapunov function $V$ are defined as
\begin{align}
  V_o &= r\epsilon^T P_o \epsilon
  \ \ ; \ \ 
        V_c = \frac{1}{2}x_1^2 + r\eta^T P_c\eta.
        \ \ ; \ \ 
V =cV_o + V_c
\label{Vdefn}
\end{align}
where $c$ is a positive constant to be picked further below. From \eqref{sfobs}, \eqref{sfetamatdyn_of}, and \eqref{Vdefn}, we have
\begin{align}
 \dot V
  &=
cr^2\epsilon^T[P_oA_o+A_o^TP_o]\epsilon  
    + 2rc\epsilon^TP_o\overline\Phi
    + x_1[\phi_1+(r\eta_2 -\zeta_1x_1-r\epsilon_2)\phi_{(1,2)}]
  \nonumber\\&\quad 
                 +r^2\eta^T[P_c A_c +A_c^T P_c]\eta
  +2r\eta^T P_c(\Phi+H[\eta_2-\epsilon_2]+\Xi)
  -2r^2\eta^T P_cG\epsilon_2
\nonumber\\&\quad 
  -\dot r c\epsilon^T[P_o\tilde D_o\!+\!\tilde D_oP_o]\epsilon
  -\dot r\eta^T[P_c \tilde D_c\!+\!\tilde D_c P_c]\eta
  \label{Vdot1_of}
\end{align}
where
$\tilde D_o=D_o- \frac{1}{2}I_{n-1}$ and
$\tilde D_c=D_c- \frac{1}{2}I_{n-1}$.

\subsubsection{Coupled Lyapunov Inequalities}
\label{sec:coupled_lyap_of}
As discussed in Section~\ref{sec:coupled_lyap_sf}, the sets of inequalities \eqref{eq:cascading_dominance_controller} and \eqref{eq:cascading_dominance_observer} are the cascading dominance conditions in the controller and observer contexts, respectively, which require the upper diagonal terms closer to the input or output, respectively, to be bigger (in a nonlinear function sense). Analogous to the design of the symmetric positive-definite matrix $P_c$ and controller gains $k_2,\ldots,k_n$ in Section~\ref{sec:coupled_lyap_sf}, the symmetric positive-definite matrix $P_o$ and observer gains $g_2,\ldots,g_n$ (which appear in the definition of the matrix $A_o$) can be designed using the constructive procedure in \cite{KK04f,KK06}. Specifically, under Assumption A1 and condition \eqref{eq:cascading_dominance_observer} in Assumption A3, $P_o,g_2,\ldots,g_n$ can be constructed (\cite{KK04f,KK06}) such that the following coupled Lyapunov inequalities are satisfied (for all $x_1\in{\cal R}$) with some positive constants $\nu_{o}$, $\tilde\nu_o$, $\underline\nu_{o}$,
and $\overline\nu_{o}$:
\begin{align}
  P_o A_o+A_o^T P_o \leq -\nu_o  I - \tilde\nu_o \phi_{(2,3)}C^TC
  \ \ \ ; \ \ 
  \underline \nu_o I \leq P_o \tilde D_o + \tilde D_o P_o\leq \overline\nu_o I
  \label{coupled_lyap2}
\end{align}
where $C=[1,0,\ldots,0]$. From Theorem 2 in \cite{KK06}, $g_2,\ldots,g_n$ can be chosen to be linear constant-coefficient combinations of $\phi_{(2,3)},\ldots,\phi_{(n-1,n)}$. Hence, using Assumption A3, a positive constant $\overline G$ can be found such that
\begin{align}
  (\sum_{i=2}^n g_i^2)^{\frac{1}{2}} \leq \overline G \phi_{(2,3)}.
  \label{Gbardefn}
\end{align}

\subsubsection{Design Freedoms}
\label{sec:design_freedoms_of}
The design freedoms appearing in the Lyapunov inequality \eqref{Vdot1_of} are the constant $c$, the function $\zeta_1$, and the dynamics of $r$. The properties to be satisfied by the choice of these design freedoms are outlined below. As in the state-feedback case in Section~\ref{sec:design_freedoms_sf}, the dynamics of the high-gain scaling parameter $r$ are picked in the following form to ensure the property that the derivative $\dot r$ is ``large'' until $r$ itself becomes ``large,'' with ``large''-ness in both cases being defined in terms of nonlinear functions constructed based on upper bounds for various terms that appear in Lyapunov inequalities written for the closed-loop system. For this purpose, the dynamics of $r$ are designed to be of the form shown in \eqref{rdot} as in the strict-feedback case; the criteria for choices of the positive quantity $a$ (a constant or a positive function of $x_1$ lower bounded by a positive constant) and the function $\Omega$ are discussed later in this section.
With the dynamics of $r$ designed as in \eqref{rdot}, \eqref{Vdot1_of} yields
\begin{align}
  \dot V
  &=
    cr^2\epsilon^T[P_oA_o+A_o^TP_o]\epsilon  
    +r^2\eta^T[P_c A_c +A_c^T P_c]\eta
    \nonumber\\&\quad 
    +a r^2 c\epsilon^T[P_o\tilde D_o+\tilde D_oP_o]\epsilon
    +a r^2\eta^T[P_c \tilde D_c+\tilde D_c P_c]\eta
    -2r^2\eta^T P_cG\epsilon_2
    \nonumber\\&\quad 
  + 2rc\epsilon^TP_o\overline\Phi
  + rx_1(\eta_2 -\epsilon_2)\phi_{(1,2)}
  +2r\eta^T P_c(\Phi+H[\eta_2-\epsilon_2]+\Xi)
  \nonumber\\&\quad
  -r(\Omega + a) c\epsilon^T[P_o\tilde D_o+\tilde D_oP_o]\epsilon
  -r(\Omega + a)\eta^T[P_c \tilde D_c+\tilde D_c P_c]\eta
  \nonumber\\&\quad 
  -\zeta_1 x_1^2 + x_1\phi_1
  \label{Vdot2_of}
\end{align}
Using \eqref{coupled_lyap} and \eqref{coupled_lyap2}, \eqref{Vdot2_of} yields
\begin{align}
  \dot V
  &\leq
    -cr^2(\nu_o-a\overline\nu_o) |\epsilon|^2 -cr^2\tilde\nu_o\phi_{(2,3)}\epsilon_2^2
    -r^2(\nu_c\phi_{(2,3)}-a\overline\nu_c)|\eta|^2
  -2r^2\eta^T P_cG\epsilon_2
  \nonumber\\&\quad 
  + 2rc\epsilon^TP_o\overline\Phi
  + rx_1(\eta_2 -\epsilon_2)\phi_{(1,2)}
  +2r\eta^T P_c(\Phi+H[\eta_2-\epsilon_2]+\Xi)
  \nonumber\\&\quad
  -r(\Omega + a) c\epsilon^T[P_o\tilde D_o+\tilde D_oP_o]\epsilon
  -r(\Omega + a)\eta^T[P_c \tilde D_c+\tilde D_c P_c]\eta
  \nonumber\\&\quad 
  -\zeta_1 x_1^2 + x_1\phi_1
  \label{Vdot3_of}
\end{align}

The terms in \eqref{Vdot2_of} and \eqref{Vdot3_of} can be considered in terms of three parts:
    (a) terms multiplied by $r^2$, i.e., the first two lines of \eqref{Vdot2_of} or the first line of \eqref{Vdot3_of}; (b)
    terms multiplied by $r$, i.e., the third and fourth lines of \eqref{Vdot2_of} or the second and third lines of \eqref{Vdot3_of}; (c)
    terms not explicitly multiplied by $r$, i.e., the last lines of \eqref{Vdot2_of} or \eqref{Vdot3_of}.

Based on \eqref{Vdot3_sf}, the primary design considerations are summarized below. These design considerations are seen to be analogous to the state-feedback case in Section~\ref{sec:design_freedoms_sf} except for the additional constant $c$ and the various additional terms in \eqref{Vdot2_of} and \eqref{Vdot3_of} that need to be addressed in the choices of $a$, $\zeta_1$, and $\Omega$.
\begin{itemize}
\item The positive quantity $a$ appearing in the dynamics of $r$ in \eqref{rdot} must be picked small enough to retain the negativeness of the terms involving quadratics of $\epsilon$ and $\eta$ and multiplied by $r^2$ in the  first line of \eqref{Vdot3_of}. For example, 
$a$ can be picked\footnote{The notations $\max(a_1,\ldots,a_n)$ and $\min(a_1,\ldots,a_n)$ indicate the largest and smallest values, respectively, among numbers $a_1,\ldots,a_n$.} such that
  \begin{align}
    a &\leq \frac{1}{2}\min(\frac{\nu_o}{\overline\nu_o},\frac{\nu_c\sigma}{\overline\nu_c}).
        \label{a_choice1}
  \end{align}
\item The positive constant $c$ must be picked large enough so that the sign-indefinite quantity (sign-indefinite since $\eta$ and $\epsilon_2$ can have any signs) appearing as the last term in the first line of \eqref{Vdot3_of} can be dominated by the other negative terms in the first line. Writing, for example, the inequality
  \begin{align}
    -2r^2\eta^TP_c G \epsilon_2 &\leq 
                                  \frac{8}{\nu_c}\phi_{(2,3)}r^2 \lambda_{max}^2(P_c)\overline G^2 \epsilon_2^2
                                  + r^2\frac{\nu_c}{8}\phi_{(2,3)}|\eta|^2,
                                  \label{c_choice_cond1}
  \end{align}
  $c$ can be picked based on this design consideration to be such that
\begin{align}
  c &\geq \frac{32 \lambda_{max}^2(P_c)\overline G^2}{3\tilde\nu_o\nu_c}.
      \label{c_choice1}
\end{align}
\item The functions $\zeta_1$ and $\Omega_1$ are to be picked such that the sign-indefinite terms in the second line and the second term in the last line of \eqref{Vdot3_of} are dominated by the terms in the third line and the first term in the last line.
\end{itemize}
The designs of $\zeta_1$ and $\Omega$ based on the design consideration above entail computing upper bounds on the terms in the second line of \eqref{Vdot3_of}. For example,
  using the bounds on the uncertain terms $\phi_i$ in Assumption~A2 and the property that $r\geq 1$, the functions $\overline\Phi_i$ defined in \eqref{Phiidefn} can be bounded as
  \begin{align}
    \!\!\!\!\!\!\!\!|\overline\Phi_i| &\leq \frac{1}{r^{i-1}}[\phi_{(i,1)}(x_1)|x_1|+\phi_{(i,2)}(x_1)|\zeta|]
                        \nonumber\\
                      &\!\!\!\!\!\!\!\!\!\!
                        + |g_i(x_1)x_1|\frac{\phi_{(1,1)}(x_1)}{\phi_{(1,2)}(x_1)}
                        + \sum_{j=2}^i\phi_{(i,j)}(x_1)[|\eta_j|+|\epsilon_j|]
                        \label{eq:Phiibount}
  \end{align}
  from which an upper bound for $|\overline\Phi|$ can be computed as
  \begin{align}
    |\overline\Phi| &\leq \frac{|x_1|}{r}[|\tilde\phi_1| + |\zeta_1(x_1)||\tilde\phi_2|]
                      + ||\tilde A(x_1)|| (|\eta| + |\epsilon|)
                      \nonumber\\
                    &\quad                    + |x_1|\frac{\phi_{(1,1)}(x_1)}{\phi_{(1,2)}(x_1)}\overline G \phi_{(2,3)}(x_1)
                      \label{eq:Phibarbound_of}
  \end{align}
  where $\tilde\phi_1$, $\tilde\phi_2$, and $\tilde A$ are defined as in Section~\ref{sec:design_freedoms_sf}. 
From \eqref{eq:Phibarbound_of}, an upper bound on
$2rc\epsilon^T P_o\overline\Phi$ can be computed. Similarly, a term such as $2r\eta^T P_cH(\eta_2-\epsilon_2)$ can be upper bounded as 
  \begin{align}
    2r\eta^T P_cH(\eta_2-\epsilon_2) &\leq 3r \lambda_{max}(P_c) \phi_{(1,2)}|\zeta_1'x_1+\zeta_1|
      [|\eta|^2+|\epsilon|^2].
      \label{bnd_ineq1}
  \end{align}
  Bounds for other terms in the second line of \eqref{Vdot3_of} can be derived similarly. Using these upper bounds, the functions $\zeta_1$ and $\Omega$ can be designed so that the terms in the third line and first term of the last line of \eqref{Vdot3_of} dominate these computed upper bounds.  While the construction above of the design freedoms $a$, $c$, $\zeta_1$, and $\Omega$ suffices for the purpose of proving asymptotic stability of the closed-loop system, the constructions of the design freedoms are conservative since they effectively capture worst-case bounds that do not take into account the specific structure of the state appearance in the uncertain terms in the system. Furthermore, the explicit computations of the upper bounds for the construction of these design freedoms entails by-hand algebraic complexity as discussed above. In comparison, it will be seen in Section~\ref{sec:design_of} that a matrix pencil based approach enables computation of the design freedoms to obtain less conservative values and with less algebraic complexity for computations.

\section{Matrix Pencil Based Controller Design Under State-Feedback Case}
\label{sec:design_sf}
In this section, we develop a matrix pencil based formulation for picking the design freedoms $a$, $\zeta_1$, and $\Omega$. The overall desired objective that we want to attain in the choice of the design freedoms is to ensure an inequality of the following form
where $V$ is the Lyapunov function defined in \eqref{Vdefn_sf}
\begin{align}
  \dot V &\leq -\kappa V
  \label{eq:Vdot}
\end{align}
with $\kappa$ being a positive constant (or a positive function of $x_1$ lower bounded by a positive constant). For this purpose, the design procedure  is structured as three steps in which each step addresses the choice of one or more of the design freedoms to achieve specific inequalities that will thereafter be formulated in terms of matrix inequalities and related matrix pencil structures and their generalized eigenvalues. In the first step, $a$ is picked based on the design consideration for $a$ as outlined in Section~\ref{sec:design_freedoms_sf} to ensure that the following inequality is satisfied with $\underline c_1$ being a constant that can be picked in the interval $(0,1)$:
\begin{align}
  0 &\geq   (1 - \underline c_1)r^2\eta^T[P_c A_c +A_c^T P_c]\eta
    +a r^2\eta^T[P_c \tilde D_c\!+\!\tilde D_c P_c]\eta.
      \label{a_cond_sf}
\end{align}
Since $P_c A_c +A_c^T P_c$ is a negative definite matrix by the construction of $P_c$, it is seen that achieving \eqref{a_cond_sf} essentially entails picking $a$ to be small enough. This is in line with the design consideration for $a$ as discussed in Section~\ref{sec:design_freedoms_sf}.
In the second step, the functions $\zeta_1$ and $\Omega$ are picked such that the following inequality is satisfied with $\underline c_2$ picked in the interval $(0,\underline c_1)$:
\begin{align}
  0 &\geq  (\underline c_1\!-\!\underline c_2)r^2 \eta^T[P_c A_c \!+\!A_c^T P_c]\eta
      + rx_1\eta_2\phi_{(1,2)}
      \nonumber\\&\quad
  +2r\eta^T P_c(\Phi+H\eta_2+\Xi)
  -r(\Omega + a) \eta^T[P_c \tilde D_c\!+\!\tilde D_c P_c]\eta
  \nonumber\\&\quad
  -(1\!-\!\underline c_2)\zeta_1 x_1^2\phi_{(1,2)} \!+\! x_1\phi_1.
  \label{zeta1_Omega_kappa_cond_sf}
\end{align}
In the third step, $\kappa$ is picked such that
\begin{align}
  0 &\geq 
  \underline c_2 r^2 \eta^T[P_c A_c \!+\!A_c^T P_c]\eta
      - \underline c_2\zeta_1 x_1^2\phi_{(1,2)}
  +\! \kappa \Big\{\frac{x_1^2}{2} \!+\! r\eta^T P_c \eta \Big\}.
  \label{eq:kappa_cond_sf}
\end{align}
Each of these steps is discussed in detail below. The design philosophy is motivated by the observation that since each of the design freedoms is a scalar and $a$ and $\Omega$ appear linearly\footnote{It is to be noted that $\zeta_1$ does not appear completely linearly in the Lyapunov inequalities since $H$ and $\Xi$ also depend on $\zeta_1$ and its derivative $\zeta_1'$; however, it will be seen that this implicit nonlinear appearance of $\zeta_1$ can be handled in an iterative approach within a matrix pencil based design formulation.}  in the Lyapunov inequalities above, the conditions above can be written in terms of quadratic forms involving some of the state variables and the choice of the design freedoms can be expressed in terms of matrix pencil based subproblems, specifically in terms of generalized eigenvalues of appropriately defined matrix pencil structures\footnote{Given square matrices $A_1$ and $A_2$, the generalized eigenvalues of the matrix pencil $A_1 - sA_2$ are defined as the values of $s$ that make $\mbox{det}(A_1-sA_2)=0$ where $\mbox{det}$ denotes the matrix determinant. The set of generalized eigenvalues of the matrix pencil $A_1-sA_2$ are denoted as $\sigma(A_1,A_2)$. The subset of these eigenvalues that are finite in magnitude are denoted as $\sigma_f(A_1,A_2)$. 
  It can be seen that when $A_1$ and $A_2$ are symmetric and at least one of the two matrices $A_1$ and $A_2$ is positive-definite (or negative-definite), the generalized eigenvalues are real numbers. Given symmetric matrices $A_1$ and $A_2$ with at least one of $A_1$ and $A_2$ being positive-definite (or negative-definite),
the minimum eigenvalue and the largest finite eigenvalue are denoted as $\sigma_{\min}(A_1,A_2)\eqdef \min(\sigma(A_1,A_2))$ and $\sigma_{\max,f}(A_1,A_2)\eqdef \max(\sigma_f(A_1,A_2))$.}.

Since $a$ appears as part of a stabilizing negative term in $\dot r$, larger values of $a$ would be desirable to reduce conservatism of the design. However, $a$ needs to be picked to be ``small enough'' as discussed as part of the first step above. Hence, the appropriate choice of $a$ entails finding the largest possible ``small enough'' value for $a$.  On the other hand, it is seen that larger values of $\zeta_1$ and $\Omega$ tend to increase the magnitudes of the overall control input signal; hence, smaller values of $\zeta_1$ and $\Omega$ would be desirable to reduce conservatism. However, as discussed as part of the second step above, it is seen that $\zeta_1$ and $\Omega$ need to be picked to be ``large enough''. Hence the suitable choices of $\zeta_1$ and $\Omega$ entails finding smallest possible ``large enough'' values for $\zeta_1$ and $\Omega$.  

\subsection{Design of $a$}
Noting that the right hand side of \eqref{a_cond_sf} is homogeneous in $r^2$ where $r\geq 1$ and writing it as a quadratic form in terms of $\eta$,  \eqref{a_cond_sf} can be written as
\begin{align}
  a Q_{a1} + Q_{a2} &\leq 0
                                                               \label{a_cond2_sf}
\end{align}
where the matrices $Q_{a1}$ and $Q_{a2}$ are given by
\begin{align}
  Q_{a1} &=
           P_c\tilde D_c + \tilde D_c P_c
  \ \ ; \ \ 
  Q_{a2} = 
           (1 - \underline c_1)(P_cA_c + A_c^T P_c).
           \label{Qa1_Qa2_sf}
\end{align}
Note that both $Q_{a1}$ and $Q_{a2}$ are known functions of $x_1$ given a choice of $\underline c_1$ in the interval $(0,1)$.
While the nominal construction of $a$ in \eqref{a_choice1_sf} is as a positive constant, $a$ can be allowed to be a function of $x_1$ as long as it has a positive lower bound, i.e., $a(x_1)\geq \underline a$ with some positive $\underline a$. Indeed, allowing $a$ to be dependent on $x_1$ can enable larger values of $a$, which is beneficial since the term involving $a$ essentially provides a stabilizing effect on the dynamics of $r$. Also, as can be seen from the discussion following equation \eqref{rdot}, a larger value of $a$ would tend to make the eventual (steady state) value of $r$ smaller, which also tends to reduce effective control gains since the definition of the control input $u$ in \eqref{sfcontrol} involves $r^n$, thus making the steady state value of $r$ smaller preferable (also, note that $r$ is monotonically non-decreasing due to the dynamics of $r$ given in \eqref{rdot}).

From the construction in Section~\ref{sec:coupled_lyap_sf}, note that $Q_{a1}$ defined in \eqref{Qa1_Qa2_sf} is a symmetric positive-definite matrix and $Q_{a2}$ is negative-definite.  Hence, analogous to the reasoning in the choice of $a$ in \eqref{a_choice1_sf} in Section~\ref{sec:design_freedoms_sf}, it is seen that picking $a$ small enough will make $aQ_{a1} + Q_{a2}$ negative semidefinite.
Hence, noting the matrix pencil structure of $aQ_{a1} + Q_{a2}$, it is seen that $a$ can simply be picked to be a function of $x_1$ obtained as
\begin{align}
  a(x_1) &= \sigma_{\min}(Q_{a2}(x_1),-Q_{a1}(x_1)).
           \label{a_choice_sf}
\end{align}
Note that none of the generalized eigenvalues of the matrix pencil $Q_{a2}-sQ_{a1}$ can be negative since with any $s<0$, the negative definiteness of $Q_{a2}-sQ_{a1}$ follows from negative definiteness of $Q_{a2}$ and positive definiteness of $Q_{a1}$. Hence, all generalized eigenvalues of the matrix pencil $Q_{a2}-sQ_{a1}$ are positive and \eqref{a_choice_sf} prescribes $a$ as the smallest of the generalized eigenvalues. Since we know that a small enough $a$ definitely exists such that $aQ_{a1}+Q_{a2}$ is negative semidefinite for all $a$ smaller than this value, it follows that $a$ computed in \eqref{a_choice_sf} is the (largest such) small enough value. Indeed, it is seen that for any $a$ smaller than the value constructed in \eqref{a_choice_sf}, it can be seen that $aQ_{a1}+Q_{a2}$ definitely satisfies \eqref{a_cond2_sf}.

\subsection{Designs of $\zeta_1$ and $\Omega$}
From \eqref{zeta1_Omega_kappa_cond_sf}, the basic strategy is to pick $\zeta_1$ and $\Omega$ large enough to make \eqref{zeta1_Omega_kappa_cond_sf} hold. However, there are several subtle details to consider. Firstly, in terms such as $2r\eta^T P_c\Phi$, the sign of elements of $\Phi$ is not known (since the only information on uncertain terms $\phi_i$ is the structure of upper bounds in Assumption~A2). Hence, when converting \eqref{zeta1_Omega_kappa_cond_sf} into a matrix pencil structure, we will need to address the uncertainties in signs of the constituent terms. One way to address this is to write quadratic forms in terms of vectors of element-wise magnitudes\footnote{Given a vector $a=[a_1,\ldots,a_m]^T$, the notation $|a|_e$ is used to denote the vector of element-wise magnitudes $[|a_1|,\ldots,|a_m|]^T$.
  Given two vectors $a$ and $b$ of the same dimension $m\times 1$, the relation $a \leq_e b$ indicates the set of element-wise inequalities $|a_i|\leq |b_i|, i=1,\ldots,m$ where $a_i$ and $b_i$ indicate the $i^{th}$ element of $a$ and $b$, respectively. Also, given a matrix $M$, the notation $|M|_e$ is used to denote the matrix of the same dimensions with each element being the magnitude of the corresponding element of $M$.} $|\eta|_e$ rather than in terms of $\eta$ itself. However, this then requires that the sign-definite terms in the first and second lines of \eqref{zeta1_Omega_kappa_cond_sf} should also be written in terms of $|\eta|_e$ rather than in terms of $\eta$ itself. For this purpose, while conservative estimates can be written such as, for example, $\eta^T(P_cA_c+A_c^TP_c)\eta\leq -\nu_c\phi_{(2,3)} |\eta|_e^T |\eta|_e$ using the coupled Lyapunov inequalities in \eqref{sec:coupled_lyap_sf}, relatively non-conservative estimates can be obtained by considering the diagonal elements of matrices $P_cA_c+A_c^TP_c$, etc. Specifically, defining\footnote{Given a matrix $M$, the notation $\diag(M)$ is used to denote the diagonal matrix whose diagonal elements are equal to the diagonal elements of $M$.}
\begin{align}
  \overline A_c &= \mbox{diag}(P_cA_c+A_c^T P_c),
\end{align}
we have $\eta^T \overline A_c\eta = |\eta|_e^T \overline A_c|\eta|_e$ since $\overline A_c$ is a diagonal matrix and considering the matrix pencil of form $P_cA_c+A_c^TP_c - s\overline A_c$, we see that defining
\begin{align}
  \delta_{A_c} &= \sigma_{\min}(P_cA_c+A_c^TP_c , \overline A_c),
\end{align}
we have the inequality
  $P_cA_c+A_c^TP_c \leq \delta_{A_c}\overline A_c$.
Similarly, defining
  $\overline D_c = \mbox{diag}(P_c\tilde D_c+\tilde D_c P_c)$
and
  $\delta_{D_c} = \sigma_{\min}(P_c\tilde D_c+\tilde D_c^TP_c , \overline D_c)$
we have the inequality
$P_c\tilde D_c+\tilde D_c^TP_c \geq \delta_{D_c}\overline D_c$.
Hence, noting that $\overline A_c$ and $\overline D_c$
are diagonal matrices, the inequality \eqref{zeta1_Omega_kappa_cond_sf} reduces to
\begin{align}
  0 &\geq  (\underline c_1\!-\!\underline c_2)r^2\delta_{A_c}|\eta|_e^T\overline A_c|\eta|_e
      + rx_1\eta_2\phi_{(1,2)}
      \nonumber\\&\quad
  +2r\eta^T P_c(\Phi+H\eta_2+\Xi)
  -r(\Omega + a) \delta_{D_c}|\eta|_e^T\overline D_c|\eta|_e
  \nonumber\\&\quad
  -(1\!-\!\underline c_2)\zeta_1 x_1^2\phi_{(1,2)} \!+\! x_1\phi_1.
  \label{zeta1_Omega_kappa_cond2_sf}
\end{align}
Note that $\overline A_c$ is negative-definite while $\overline D_c$
is positive-definite.

The second detail to consider is that different groups of terms in \eqref{zeta1_Omega_kappa_cond2_sf} involve the scaling parameter $r$ to different powers (i.e., terms multiplied by $r^2$, terms multiplied by $r$, and terms not multiplied by a power of $r$). Also, terms such as $rx_1\eta_2\phi_{(1,2)}$ cannot (in an $r$-independent way) be dominated by the negative terms $-\zeta_1 x_1^2$ and $-r(\Omega + a)\delta_{D_c}|\eta|_e^T\overline D_c|\eta|_e$ even if $\zeta_1$ and $\Omega$ are picked large. Note that both $\zeta_1$ and $\Omega_1$ are required to be functions of $x_1$ and cannot depend on $r$. To address these considerations, the inequality in \eqref{zeta1_Omega_kappa_cond2_sf} is formulated as a quadratic form in terms of the expanded vector
\begin{align}
  X &= [|x_1|,\sqrt{r}|\eta|_e^T,r|\eta|_e^T]^T.
\end{align}
Next, note that $\zeta_1$ appears within the definitions of $H$ and $\Xi$ and that therefore, $\zeta_1$ effectively appears nonlinearly in \eqref{zeta1_Omega_kappa_cond2_sf}. However, from \cite{KK04f}, it is known that a $\zeta_1$ of the form (with $c_{\zeta_1}$ being any positive constant)
\begin{align}
  \!\!\!\zeta_1(x_1) &\!=\! a_{\zeta_1} \Big[1 \!+\! c_{\zeta_1}\Gamma(x_1)\Big]
  \label{eq:zeta1form_sf}
\end{align} 
is admissible with some large enough positive constant $a_{\zeta_1}$ where $\Gamma(x_1) = \sum_{i=2}^n\sum_{j=1}^i \phi_{(i,j)}(x_1)$. 
To determine the appropriate value of $a_{\zeta_1}$, a feasibility test based on the computed value of $\Omega$ is used as discussed below.
With $\zeta_1$ picked to be of the form \eqref{eq:zeta1form_sf}, it can be seen that the inequality \eqref{zeta1_Omega_kappa_cond2_sf} can be written using quadratic forms in terms of $X$ and using the fact that $r\geq 1$ as
\begin{align}
  \Omega Q_{\Omega 1} + Q_{\Omega 2} &\leq 0
                                                 \label{zeta1_Omega_kappa_cond3_sf}
\end{align}
where $Q_{\Omega 1} = \diag(0, -\delta_{D_c}\overline D_c, 0)$
with the $0$ elements in the block diagonal matrix $Q_{\Omega 1}$ being defined to be of dimensions to be compatible with the quadratic form structure in terms of $X$. 
$Q_{\Omega 2}$ is a matrix function of $x_1$ of dimension $[2(n-1)+1]\times [2(n-1)+1]$ that can be obtained using the following inequalities that can be written using Assumption~A2 and \eqref{Phidefn}--\eqref{Xi_defn}:
\begin{align}
  rx_1\eta_2\phi_{(1,2)} &\leq r|x_1|B_1|\eta|_e\phi_{(1,2)}
                           \label{Q_omega2_bnd1} \\
  2r\eta^TP_c\Phi &\leq 2|\eta|_e^T|P_c|_e [|\tilde\phi_1|_e+|\tilde\phi_2|_e|\zeta_1|]|x_1|+2r|\eta|_e^T|P_c|_e\tilde A |\eta|_e
                    \label{Q_omega2_bnd2} \\
  2r\eta^TP_cH\eta_2 &\leq 2r |\eta|_e^T |P_cHB_1^T|_e|\eta|_e
                    \ \ ; \ \ 
  2r\eta^TP_c\Xi \leq 2 |\eta|_e^T |P_c\overline\Xi B_1|_e|x_1|
                   \label{Q_omega2_bnd3}
\end{align}
where
\begin{align}
  \overline\Xi &= (\phi_{(1,1)}+|\zeta_1|\phi_{(1,2)})|\zeta_1'x_1+\zeta_1|,
\end{align}
$\tilde\phi_1$ and $\tilde\phi_2$ are as given in \eqref{tilde_phi1_phi2_defn},
and
$B_1$ is the $(n-1)\times 1$ vector given by $[1,0,\ldots,0]^T$. Note that $B_1^T\eta=\eta_2$.
Using \eqref{bnd_ineq2_sf}, \eqref{zeta1_Omega_kappa_cond2_sf}, and \eqref{Q_omega2_bnd1}--\eqref{Q_omega2_bnd3}, $Q_{\Omega 2}$ can be written as
\begin{align}
  Q_{\Omega 2} &= \left[\begin{array}{ccc}
                          0 & |\overline\Xi B_1^TP_c|_e + \tilde\phi_{\zeta}^T|P_c|_e  & \frac{1}{2}B_1^T\phi_{(1,2)} \\
                          |\overline\Xi P_cB_1|_e+ |P_c|_e \tilde\phi_{\zeta} & 0& 0 \\
                          \frac{1}{2}B_1\phi_{(1,2)} & 0 & 0
                        \end{array}
                                                           \right]
                                                           \nonumber\\
                                                           &\quad +
                                                             \diag\Big(-(1-\underline c_2)\zeta_1\phi_{(1,2)}+\phi_{(1,1)},
                                                             \nonumber\\
               &\ \ \qquad\qquad
                 -a\delta_{D_c}\overline D_c+ |P_cHB_1^T|_e+|B_1H^TP_c|_e+|P_c|_e\tilde A+\tilde A^T |P_c|_e ,
                 \nonumber\\
               &\ \ \qquad\qquad
                 (\underline c_1-\underline c_2)\delta_{A_c}\overline A_c
                 \Big).
                 \label{Q_omega2_sf}
\end{align}
where
  $\tilde\phi_{\zeta} = [|\tilde\phi_1|_e+|\tilde\phi_2|_e|\zeta_1|]$.
Noting the structure of the matrix pencil $\Omega Q_{\Omega 1} + Q_{\Omega 2}$, $\Omega$ is designed as
\begin{align}
  \Omega &= \sigma_{max,f}(Q_{\Omega 2}, -Q_{\Omega 1}).
  \label{eq:Omega_sf}
\end{align}
The non-negativeness of $\Omega$ forms a simple test as to whether $a_{\zeta_1}$ is large enough. Specifically, if the obtained $\Omega$ from \eqref{eq:Omega_sf} is not non-negative, the $a_{\zeta_1}$ is increased (e.g., by a constant factor such as $1+q_a$ with some small $q_a > 0$) until a non-negative $\Omega$ is obtained. Since $a_{\zeta 1}$ and $\Omega$ are computed at run-time, the recomputation of $a_{\zeta 1}$ and $\Omega$ can result in non-smoothness of the solution trajectories. However, such points of non-smoothness occur at most a finite number of times. To see this, note that it is known from the reasoning in Section~\ref{sec:design_freedoms_sf} and the detailed construction in \cite{KK04f} that a value of $a_{\zeta_1}$ 
that is large enough to make a choice of $\Omega$ feasible definitely exists. Hence, whenever this iterative procedure for computation of $a_{\zeta 1}$ and $\Omega$ is activated due to a negative $\Omega$, the iterative procedure definitely terminates within a finite number of increases of $a_{\zeta 1}$. Since the closed-loop dynamics equations will have discontinuities at points of time at which the choice of $a_{\zeta 1}$ is revised upward, the solutions of the differential equations are understood in the sense of Filippov~\cite{filippov} Furthermore, by continuity of underlying functions, this iterative procedure will be activated at most an explicitly bounded number of times, implying that non-smoothness points in the solution trajectories occur only at most an explicitly bounded number of times and solutions to the closed-loop system are defined by simple concatenation.

\subsection{Design of $\kappa$}
Once the design freedoms $a$, $\zeta_1$, and $\Omega$ are picked as discussed above, $\kappa$ can be easily found to satisfy inequality \eqref{eq:kappa_cond_sf}. The small enough $\kappa$ for this purpose is found from the matrix pencil structure of \eqref{eq:kappa_cond_sf} by writing \eqref{eq:kappa_cond_sf} as a quadratic form in terms of $[x_1,\sqrt{r}\eta^T]^T$; noting that $r\geq 1$ while $P_cA_c+A_c^TP_c$ is negative-definite, we obtain the matrix inequality
\begin{align}
  \kappa Q_{\kappa 1}+Q_{\kappa 2}&\leq 0
                                    \label{eq:kappa_cond2}
\end{align}
where
\begin{align}
  Q_{\kappa 1} &= \left[\begin{array}{cc}
                          \frac{1}{2} & 0 \\
                          0 & P_c
                        \end{array}
                 \right]\ \ ; \ \ 
  Q_{\kappa 2} = \left[\begin{array}{cc}
                          -\underline c_2 \zeta_1\phi_{(1,2)} & 0 \\
                          0 & \underline c_2[P_cA_c+A_c^TP_c]
                        \end{array}
                              \right].
\end{align}
Note that $Q_{\kappa 1}$ is a symmetric positive-definite matrix while $Q_{\kappa 2}$ is a symmetric negative-definite matrix. Hence, \eqref{eq:kappa_cond2} can be satisfied simply by picking $\kappa$ small enough.
Hence, from \eqref{eq:kappa_cond2}, it is seen that the largest possible small enough value for $\kappa$ is given by
\begin{align}
  \kappa &= \sigma_{min}(Q_{\kappa 2}(x_1),-Q_{\kappa 1}(x_1))
  \label{eq:kappa}
\end{align}

Note that $\kappa$ is not required in the control implementation, but is only used to analyze convergence from \eqref{eq:Vdot}. Also, note that as with the choice of $a$, it is acceptable for $\kappa$ to be a function of $x_1$ as long as it is lower bounded by a positive constant since a state dependence of $\kappa$ will not create any additional terms in $\dot V$ and stability and convergence properties follow from \eqref{eq:Vdot}.

\subsection{Simulation Studies}
\label{sec:example_sf}
\begin{figure}
    \centerline{\includegraphics[width=5in]{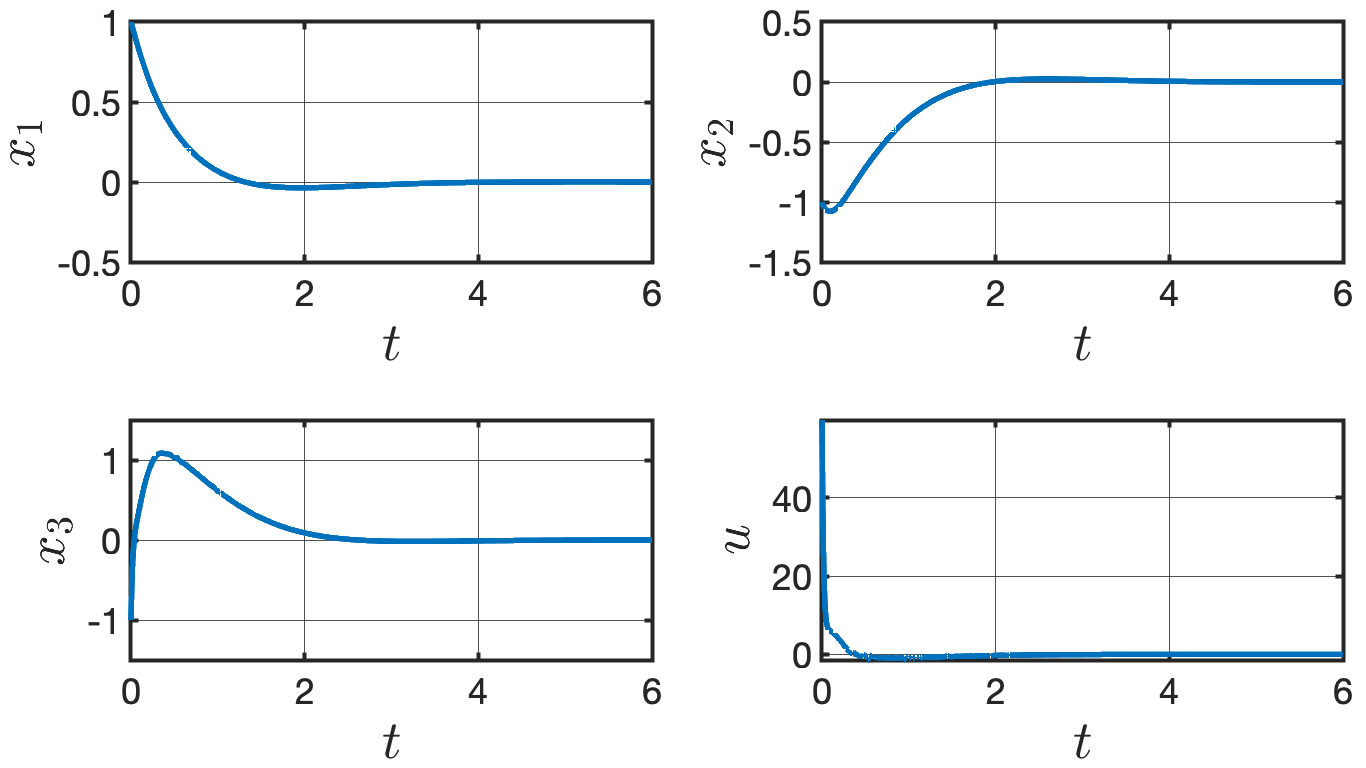}}
    \caption{Simulations for the closed-loop system with the designed state-feedback controller.}
    \label{fig:sim_sf}
  \end{figure}

To study the efficacy of the proposed approach, consider the following example system:
\begin{align}
  \dot x_1 &= (1+x_1^2)x_2 \ \ ; \ \ y = x_1\\
  \dot x_2 &= (1+x_1^4)x_3 + x_1^3 x_2 \cos(x_3)\\
  \dot x_3 &= u + x_1^2 x_2.
\end{align}
Here, $\phi_{(1,2)}(x_1)=1+x_1^2$, $\phi_{(2,3)}(x_1)=1+x_1^4$, and $\mu_0(x_1)=1$. This system is easily seen to satisfy the Assumptions~A1--A3 given in Section~\ref{sec:formulation}. In particular, $\sigma=1$,
$\phi_{(1,1)} = \phi_{(2,1)} = \phi_{(3,1)} = \phi_{(3,3)} = 0$,
$\phi_{(2,2)} = x_1^3$, and
$\phi_{(3,2)} = x_1^2$.
Using the constructive procedure in \cite{KK04f,KK06},
a symmetric positive-definite matrix $P_c$ and functions $k_2$ and $k_3$ are found
to satisfy the coupled Lyapunov inequalities \eqref{coupled_lyap} as $P_c=\tilde a_c\left[\begin{array}{cc}2 & 0.5 \\ 0.5&1 \end{array}\right]$, 
  $k_2=8\phi_{(2,3)}$, and $k_3=15\phi_{(2,3)}$,
  and with $\nu_c=1.397 \tilde a_c$, $\underline\nu_c=1.382\tilde a_c$, and $\overline\nu_c=3.618\tilde a_c$ with $\tilde a_c$ being any positive constant. 
  Picking $\tilde a_c=0.05$, the design parameters in the proposed control design approach are picked as 
  $\underline c_1 = 0.1$,
  $\underline c_2 = 0.0001$,
  $c_{\zeta_1}=0.001$, $q_a=0.02$.
  The initial value of $a_{\zeta_1}$ is picked as 1.
  For simulation studies, the initial condition for the system state vector $[x_1,x_2,x_3]^T$ is specified as $[1,-1,-1]^T$. The closed-loop trajectories and the control input signal are shown in Figure~\ref{fig:sim_sf}.  
  As a quick point of comparison with the ``macroscopic'' conservative approach in the controller construction in prior results (e.g., \cite{KK04f}), it can be seen that the construction of $a$ as in 
  \eqref{a_choice1_sf} provides $a=0.1931$ while the values of $a$ dynamically computed in the proposed approach vary between 0.895 and 1.79, yielding a far less conservative value for $a$ (note that larger values for $a$ are desirable since $a$ acts as a stabilizing term in the dynamics of $r$). A comparison of $\Omega$ shows a similar trend with the proposed approach yielding less conservative (lower) values compared to the more conservative  ``macroscopic'' approach outlined in Section~\ref{sec:design_freedoms_sf}, which is developed in more detail in \cite{KK04f}. While simulations with the design freedoms computed as in Section~\ref{sec:design_freedoms_sf} are omitted here for brevity, a comparison shows that the matrix pencil based approach provides significantly smaller overshoots in the transients and smaller control magnitudes due to less conservative values of design freedoms.

\section{Matrix Pencil Based Controller Design Under Output-Feedback Case}
\label{sec:design_of}
Following a similar approach as in the state-feedback case in Section~\ref{sec:design_sf}, we develop a matrix pencil based formulation in this section for picking the design freedoms in the dynamic output-feedback controller designed in Section~\ref{sec:scaling_control_of}. In the output-feedback case, the design freedoms are $c$, $a$, $\zeta_1$, and $\Omega$. The overall objective that is desired to be attained through  appropriate choice of the design freedoms is to ensure an inequality of the form \eqref{eq:Vdot}
with $\kappa$ being a positive constant (or a positive function of $x_1$ lower bounded by a positive constant).
The Lyapunov function $V$ is as given in \eqref{Vdefn} and includes quadratics in $x_1$, $\eta$, and $\epsilon$.
The design procedure is structured below in terms of four steps in each of which values for one or more of the design freedoms are chosen to satisfy intermediate inequalities written to eventually attain an inequality of the form \eqref{eq:Vdot}. In the first step, $c$ is picked to be a positive constant such that the following inequality is satisfied:
\begin{align}
  0
  &\geq
    (1-\underline c_1)\{cr^2\epsilon^T[P_oA_o+A_o^TP_o]\epsilon
    +r^2\eta^T[P_c A_c +A_c^T P_c]\eta\}
    - 2r^2\eta^T P_cG\epsilon_2
    \label{c_cond}
\end{align}
with $\underline c_1$ being a constant that can be picked in the interval $(0,1)$.
In the second step, $a$ is picked (to be a function of $x_1$ bounded below by a positive constant) to ensure that the following inequality is satisfied with $\underline c_2$ being a constant that can be picked in the interval $(0,\underline c_1)$:
\begin{align}
  0 &\geq   (\underline c_1 - \underline c_2)\{cr^2\epsilon^T[P_oA_o+A_o^TP_o]\epsilon
      +r^2\eta^T[P_c A_c +A_c^T P_c]\eta\}
      \nonumber\\
    &\ 
      +a r^2 c\epsilon^T[P_o\tilde D_o\!+\!\tilde D_oP_o]\epsilon
      \!+\!a r^2\eta^T[P_c \tilde D_c\!+\!\tilde D_c P_c]\eta.
      \label{a_cond}
\end{align}
In the third step, $\zeta_1$ and $\Omega$ are picked as functions of $x_1$ such that the following inequality is satisfied with $\underline c_3$ picked in the interval $(0,\underline c_2)$:
\begin{align}
  0 &\geq  (\underline c_2\!-\!\underline c_3)r^2 \{c\epsilon^T[P_oA_o\!+\!A_o^TP_o]\epsilon
      +\eta^T[P_c A_c \!+\!A_c^T P_c]\eta\}
      \nonumber\\
    &\quad
      + 2rc\epsilon^TP_o\overline\Phi
      + rx_1(\eta_2 -\epsilon_2)\phi_{(1,2)}
  +2r\eta^T P_c(\Phi+H[\eta_2-\epsilon_2]+\Xi)
  \nonumber\\&\quad
  -r(\Omega + a) \{c\epsilon^T[P_o\tilde D_o\!+\!\tilde D_oP_o]\epsilon
  +\eta^T[P_c \tilde D_c\!+\!\tilde D_c P_c]\eta\}
  \nonumber\\&\quad
  -(1\!-\!\underline c_3)\zeta_1 x_1^2\phi_{(1,2)} \!+\! x_1\phi_1.
  \label{zeta1_Omega_kappa_cond}
\end{align}
In the fourth step, $\kappa$ is picked to be a function of $x_1$ bounded below by a positive constant such that
\begin{align}
  0 &\geq 
  \underline c_3 r^2 \{c\epsilon^T[P_oA_o\!+\!A_o^TP_o]\epsilon
      +\eta^T[P_c A_c \!+\!A_c^T P_c]\eta\}
  \nonumber\\&\quad\!\!\!\!\!
  - \underline c_3\zeta_1 x_1^2\phi_{(1,2)}
  +\! \kappa \Big\{\frac{x_1^2}{2} \!+\! r\eta^T P_c \eta \!+\! r\epsilon^T P_o \epsilon\Big\}.
  \label{eq:kappa_cond}
\end{align}
The matrix pencil structures arising in each of these steps are discussed in detail below. While $a$, $\zeta_1$, $\Omega_1$, and $\kappa$ can be functions of $x_1$, $c$ is required to be a constant since any state dependence of $c$ will generate additional terms in $\dot V$ thereby modifying \eqref{Vdot1_of} and the following analysis. Specifically, since $c$ appears as a coefficient in the definition of the Lyapunov function $V$ in \eqref{Vdefn}, the derivative of $V$ will involve additional terms from $\dot c$ if $c$ is not a constant. To avoid such additional terms that would invalidate the following analysis, the choice of $c$ is constrained to be a constant below. To reduce conservatism, it is desirable that $a$ is picked as large as possible while $c$, $\zeta_1$, and $\Omega$ are picked as small as possible. However, intrinsically, the requirements on these design freedoms as seen in the inequalities above are that $a$ should be picked ``small enough'' while $c$, $\zeta_1$, and $\Omega$ are picked ``large enough''. Hence, in the matrix pencil based design below, we will seek to determine the largest possible ``small enough'' value for $a$ and smallest possible ``large enough'' values for $c$, $\zeta_1$, and $\Omega$. The parameter $\kappa$ is not required in the controller implementation but appears in the stability and convergence analysis from \eqref{eq:Vdot}. To obtain a less conservative estimate for the convergence rate of $V$, we would want to find the largest possible estimate for $\kappa$.

\subsection{Design of $c$}
Noting that the right hand side of \eqref{c_cond} is homogeneous in $r^2$ where $r\geq 1$ and writing it as a quadratic form in terms of $[\epsilon^T,\eta^T]^T$, inequality \eqref{c_cond} can be equivalently written as
\begin{align}
  0 &\geq c
      \left[\begin{array}{cc}(1-\underline c_1)(P_oA_o+A_o^T P_o) & 0 \\ 0 & 0\end{array}\right]
    +
    \left[\begin{array}{cc} 0 & -(P_c G B_1)^T \\ -P_c G B_1 & (1-\underline c_1)(P_cA_c+A_c^T P_c)\end{array}\right]
                                                               \label{c_cond2}
\end{align}
where $B_1$ denotes the $1\times(n-1)$ row vector with a 1 in the first element and zeros elsewhere.
In \eqref{c_cond2}, the parts shown as 0 denote, as per the standard notation, blocks of compatible dimensions based on the other shown parts of the matrices (i.e., the 0 blocks in \eqref{c_cond2} are all of dimension $(n-1)\times (n-1)$). The right hand side of \eqref{c_cond2} is of the form $c Q_{c1}(x_1) + Q_{c2}(x_1)$ with $Q_{c1}$ and $Q_{c2}$ being known and completely determined matrices (with $c_1$ having been chosen to be any constant in the interval $(0,1)$). Analogous to the reasoning in the choice of $c$ in \eqref{c_choice1} in Section~\ref{sec:design_freedoms_of} from \eqref{c_choice_cond1}, it is known that $c$ can be constructed (for example similar to \eqref{c_choice1}) to make \eqref{c_cond} hold. However, such a construction of $c$ is conservative since it is purely in terms of ``macroscopic'' quantities such as $\lambda_{max}(P_c)$ and $\tilde\nu_o$ and ignores the finer structure of $P_c$, $P_o$, etc. Instead, by directly addressing the requirement that we want to make $c Q_{c1}(x_1) + Q_{c2}(x_1)$ negative semidefinite as shown in \eqref{c_cond2}, a much less conservative estimate of $c$ can be found.  To remove the dependence on $x_1$, note from Assumption~A3 (cascading dominance) that the upper diagonal terms $\phi_{(i,i+1)}$ are comparable (in a nonlinear function sense) and also note from Sections~\ref{sec:coupled_lyap_sf} and \ref{sec:coupled_lyap_of} that the functions $g_2,\ldots,g_n$, and $k_2,\ldots,k_n$ can be picked to be linear constant-coefficient combinations of the upper diagonal terms 
$\phi_{(2,3)},\ldots,\phi_{(n-1,n)}$. Hence, dividing \eqref{c_cond2} throughout by $\phi_{(2,3)}$, each matrix appearing in the resulting inequality varies in a polytopic set whose vertices can be computed in terms of the constants $\overline\rho_i$ and $\underline\rho_i, i=3,\ldots,n-1$, and the coefficients in the designs of functions $g_2,\ldots,g_n$ and $k_2,\ldots,k_n$. In the special case when all $\phi_{(i,i+1)}$ are identical (except for possibly different constant coefficients), the polytope reduces to a single point. In either case, the resulting system of equations (diagonally concatenated over vertices of the polytope or from the single value when the polytope reduces to a single point) can be written in the form $0\geq c \overline Q_{c1} + \overline Q_{c2}$ with $\overline Q_{c1}$ and $\overline Q_{c2}$ being constant matrices. This is a matrix pencil with coefficient $c$. Note that $\overline Q_{c1}$ and $\overline Q_{c2}$ are both symmetric matrices (hence, eigenvalues of $c\overline Q_{c1} + \overline Q_{c2}$ are all real) and $\overline Q_{c1}$ is negative semidefinite. From \eqref{c_choice1}, we know that picking $c$ large enough will definitely satisfy this inequality. Hence, when $c\rightarrow \infty$, we know that $c\overline Q_{c1} + Q_{c2}$ tends to a negative semidefinite matrix. Therefore, noting the matrix pencil structure of $c\overline Q_{c1} + \overline Q_{c2}$, it is sufficient to pick 
\begin{align}
  c &= \sigma_{\max, f}(\overline Q_{c2},-\overline Q_{c1})
      \label{c_choice}
\end{align}
since this is, by definition, the largest finite value of $c$ for which $c\overline Q_{c1} + \overline Q_{c2}$ has an eigenvalue at 0. Note that all eigenvalues of $c\overline Q_{c1} + \overline Q_{c2}$ are real for all positive $c$. Since by the reasoning above, a large enough value of $c$ is known to exist such that for all $c$ larger than this value, $c\overline Q_{c1} + \overline Q_{c2}$ is negative semidefinite, this implies that \eqref{c_choice} must be the (smallest such) large enough value.

\subsection{Design of $a$}
Similar to the analysis in the design of $c$ above, the inequality \eqref{a_cond} which is also homogeneous in $r^2$ can be equivalently written in the form \eqref{a_cond2_sf} where
\begin{align}
  Q_{a1} &= \diag\Big(c(P_o\tilde D_o + \tilde D_o P_o), (P_c\tilde D_c + \tilde D_c P_c)\Big)
                                                                                \label{Qa1}
  \\
  Q_{a2} &= 
           (\underline c_1 - \underline c_2)\diag\Big(c(P_oA_o + A_o^T P_o), (P_cA_c + A_c^T P_c)\Big).
                                                                                                     \label{Qa2}
\end{align}
Both $Q_{a1}$ and $Q_{a2}$ are known functions of $x_1$ given a choice of $\underline c_2$ in the interval $(0,\underline c_1)$.
The $c$ designed above was intrinsically required to be a constant since it is a coefficient in the definition of the Lyapunov function $V$ in  \eqref{Vdefn} and a time or state dependence of $c$ will result in additional terms in the Lyapunov inequality \eqref{Vdot1_of} and the Lyapunov inequalities following from \eqref{Vdot1_of}. However, $a$ which is simply a coefficient appearing in the dynamics of the scaling parameter $r$ in \eqref{rdot} can very reasonably be a function of $x_1$ as long as it has a positive lower bound, i.e., $a(x_1)\geq \underline a$ with some positive $\underline a$. Such a state dependence of $a$ does not introduce any new terms into the Lyapunov inequalities and does not affect the stability analysis. On the other hand, by allowing $a$ to be dependent on $x_1$, larger values of $a$ could be possible to use, thereby facilitating the beneficial stabilizing effect in the dynamics of $r$ and making the eventual (steady state) value of $r$ smaller therefore also reducing the effective control gains and control input magnitudes. 

From the construction in Section~\ref{sec:coupled_lyap_of}, it is seen that $Q_{a1}$ is a symmetric positive-definite matrix while $Q_{a2}$ is  negative-definite.  Hence, it is seen that $aQ_{a1} + Q_{a2}$ can be ensured to be negative semidefinite by picking $a$ to be ``small enough''.
Hence, from the matrix pencil structure of $aQ_{a1} + Q_{a2}$, it is seen that $a$ can be picked to be a function of $x_1$ as \eqref{a_choice_sf} with the matrices $Q_{a1}$ and $Q_{a2}$ as defined above in \eqref{Qa1} and \eqref{Qa2}. Analogous to the reasoning in the state-feedback case, it is seen that \eqref{a_choice_sf} provides the largest ``small enough'' value for $a$ that is such that $aQ_{a1}+Q_{a2}$ is negative semidefinite for all $a$ smaller than this value.

\subsection{Designs of $\zeta_1$ and $\Omega$}
To address the uncertainty in the signs of elements of $\overline\Phi$, a similar approach is used as in the state-feedback case in Section~\ref{sec:design_sf}. The structure of the appearance of $\zeta_1$ and $\Omega$ in the inequality \eqref{zeta1_Omega_kappa_cond} indicates that $\zeta_1$ and $\Omega$ are to be picked ``large enough'' to ensure that the inequality \eqref{zeta1_Omega_kappa_cond} is satisfied.  To address the uncertainty in the signs of the elements of $\overline\Phi$, the quadratic forms for the purpose of designing $\zeta_1$ and $\Omega$ are written in terms of vectors of element-wise magnitudes $|\eta|_e$ and $|\epsilon|_e$ rather than in terms of $\eta$ and $\epsilon$ themselves. Defining
\begin{align}
  \overline A_c &= \mbox{diag}(P_cA_c+A_c^T P_c) \ \ , \ \ \overline A_o = \mbox{diag}(P_oA_o+A_o^T P_o) \\
  \overline D_o &= \mbox{diag}(P_o\tilde D_o+\tilde D_o P_o) \ \ , \ \ \overline D_c = \mbox{diag}(P_c\tilde D_c+\tilde D_c P_c) \\
  \delta_{A_o} &= \sigma_{\min}(P_oA_o+A_o^TP_o , \overline A_o) \ \ , \ \ \delta_{A_c} = \sigma_{\min}(P_cA_c+A_c^TP_c , \overline A_c) \\
  \delta_{D_o} &= \sigma_{\min}(P_o\tilde D_o+\tilde D_o^TP_o , \overline D_o) \ \ , \ \ \delta_{D_c} = \sigma_{\min}(P_c\tilde D_c+\tilde D_c^TP_c , \overline D_c),
\end{align}
we have the inequalities $P_oA_o+A_o^TP_o \leq \delta_{A_o}\overline A_o$,
$P_cA_c+A_c^TP_c \leq \delta_{A_c}\overline A_c$,
$P_o\tilde D_o+\tilde D_o^TP_o \geq \delta_{D_o}\overline D_o$, and
$P_c\tilde D_c+\tilde D_c^TP_c \geq \delta_{D_c}\overline D_c$.
Hence, noting that $\overline A_o$, $\overline A_c$, $\overline D_o$, and $\overline D_c$
are diagonal matrices, the inequality \eqref{zeta1_Omega_kappa_cond} reduces to
\begin{align}
  0 &\geq  (\underline c_2 - \underline c_3) \{cr^2\delta_{A_o}|\epsilon|_e^T\overline A_o|\epsilon|_e
      +r^2\delta_{A_c}|\eta|_e^T\overline A_c|\eta|_e \}
      \nonumber\\
    &\quad
      + 2rc\epsilon^TP_o\overline\Phi
      + rx_1(\eta_2 -\epsilon_2)\phi_{(1,2)}
  +2r\eta^T P_c(\Phi+H[\eta_2-\epsilon_2]+\Xi)
  \nonumber\\&\quad
  -r(\Omega + a)\{ c\delta_{D_o} |\epsilon|_e^T\overline D_o |\epsilon|_e
  + \delta_{D_c}|\eta|_e^T\overline D_c|\eta|_e\}
  -(1-\underline c_3)\zeta_1 x_1^2\phi_{(1,2)} + x_1\phi_1
  \label{zeta1_Omega_kappa_cond2}
\end{align}
Note that $\overline A_o$ and $\overline A_c$ are negative-definite while $\overline D_o$ and $\overline D_c$
are positive-definite.

To address the fact that different groups of terms in \eqref{zeta1_Omega_kappa_cond2} involve the dynamic high-gain scaling parameter $r$ to different powers and that terms such as $rx_1(\eta_2 -\epsilon_2)\phi_{(1,2)}$ cannot (in an $r$-independent way) be dominated by the negative terms $-\zeta_1 x_1^2$ and $-r(\Omega + a)\{ c\delta_{D_o} |\epsilon|_e^T\overline D_o |\epsilon|_e
+ \delta_{D_c}|\eta|_e^T\overline D_c|\eta|_e\}$ even if the $x_1$-dependent functions $\zeta_1$ and $\Omega$ are picked large, the inequality in \eqref{zeta1_Omega_kappa_cond2} is formulated as a quadratic form in terms of the expanded vector
\begin{align}
  X &= [|x_1|,\sqrt{r}|\epsilon|_e^T,\sqrt{r}|\eta|_e^T,r|\epsilon|_e^T,r|\eta|_e^T]^T.
\end{align}
To address the nonlinear dependence on $\zeta_1$ in \eqref{zeta1_Omega_kappa_cond2} due to the fact that $H$ and $\Xi$ depend on $\zeta_1$, a similar strategy is used as in the state-feedback case in Section~\ref{sec:design_sf}.
From \cite{KK04f}, it is known that a $\zeta_1$ of the form (with $c_{\zeta_1}$ being any positive constant)
\begin{align}
  \!\!\!\zeta_1(x_1) &\!=\! a_{\zeta_1} \Big[1 \!+\! c_{\zeta_1}\phi_{(1,2)}^2 \!+\! c_{\zeta_1}\frac{|G|^2}{\phi_{(1,2)}^2}\Gamma^2(x_1) \!+\! c_{\zeta_1}\Gamma(x_1)\Big]
  \label{eq:zeta1form}
\end{align} 
is admissible with some large enough positive constant $a_{\zeta_1}$ where $\Gamma(x_1) = \sum_{i=2}^n\sum_{j=1}^i \phi_{(i,j)}(x_1)$. Hence, a two-step iterative procedure can be used to determine appropriate choices for  $\zeta_1$ and $\Omega$. Picking $c_{\zeta 1}$ to be an arbitrary positive constant and picking $a_{\zeta 1}$ to be a small positive constant, $\Omega$ is found as discussed below. If a non-negative value for $\Omega$ is found, then evidently, the value chosen for $a_{\zeta 1}$ is large enough (although perhaps larger than strictly required). However, if a negative value for $\Omega$ is obtained, then the chosen value for $a_{\zeta 1}$ is too small. In that case, the value of $a_{\zeta 1}$ is increased by a by a constant factor $1+q_a$ with $q_a$ being a small positive constant and the computation for $\Omega$ is repeated. This iterative procedure is performed until a non-negative value for $\Omega$ is obtained. 
To compute $\Omega$ given a particular choice of $a_{\zeta_1}$, the inequality \eqref{zeta1_Omega_kappa_cond2} is written using quadratic forms in terms of $X$ and using the fact that $r\geq 1$ in the form \eqref{zeta1_Omega_kappa_cond3_sf}
with $Q_{\Omega 1} = \diag(0,-c\delta_{D_o}\overline D_o, -\delta_{D_c}\overline D_c, 0,0)$ and
$Q_{\Omega 2}$ being a matrix function of $x_1$ of dimension $[4(n-1)+1]\times [4(n-1)+1]$ that can be written using Assumption~A2 and \eqref{H_defn}, \eqref{Xi_defn}, and \eqref{G_Phidefn}.
Thereafter, $\Omega$ is designed to be of the form shown in \eqref{eq:Omega_sf}.

\subsection{Design of $\kappa$}
Once the design freedoms $c$, $a$, $\zeta_1$, and $\Omega$ are picked as discussed above, the choice of $\kappa$ is essentially to pick a small enough value such that the inequality \eqref{eq:kappa_cond} is satisfied. From the matrix pencil structure of \eqref{eq:kappa_cond} when written as a quadratic form in terms of $[x_1,\sqrt{r}\epsilon^T,\sqrt{r}\eta^T]^T$, it is seen that $\kappa$ can be picked as shown in \eqref{eq:kappa} with $Q_{\kappa 1}$ and $Q_{\kappa 2}$ given by
\begin{align}
  Q_{\kappa 1} &= \diag\Big(\frac{1}{2}, P_o, P_c\Big)
\\
  Q_{\kappa 2} &= \diag\Big(-\underline c_3\zeta_1\phi_{(1,2)}, \underline c_3c[P_oA_o+A_o^TP_o], \underline c_3[P_cA_c+A_c^TP_c] \Big).
\end{align}

\subsection{Simulation Studies}
\label{sec:example_of}
\begin{figure}
    \centerline{\includegraphics[width=5.6in]{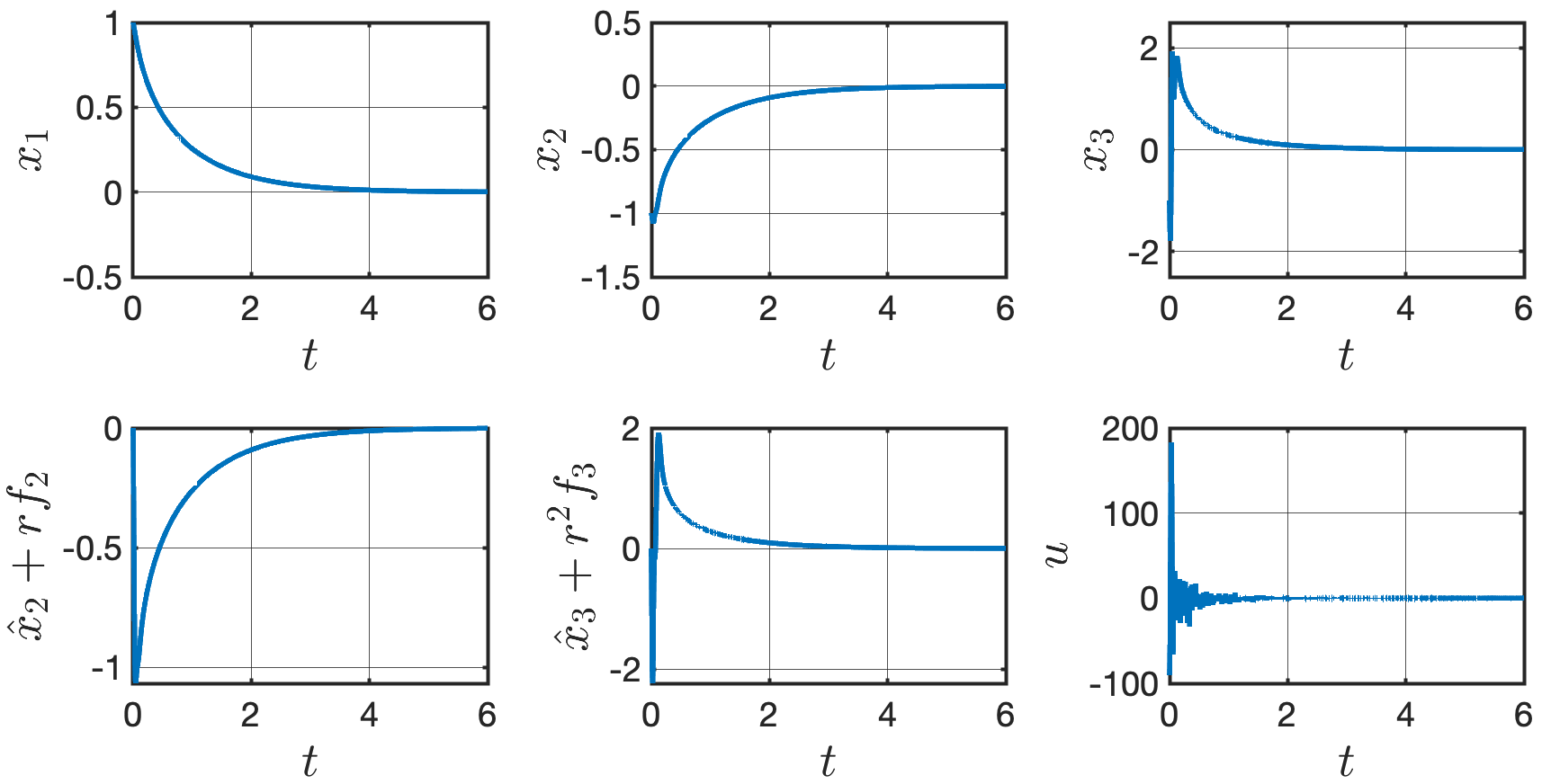}}
    \caption{Simulations for the closed-loop system with the designed output-feedback controller.}
    \label{fig:sim_of}
  \end{figure}

To study the closed-loop system behavior under the designed dynamic output-feedback controller, the same example as in Section~\ref{sec:example_sf} is considered. The constructed values of the symmetric positive matrix $P_c$, functions $k_2$ and $k_3$, and associated constants $\nu_c$, $\underline\nu_c$, $\overline\nu_c$, and $\tilde a_c$ are as in Section~\ref{sec:example_sf}.
Using the constructive procedure in \cite{KK04f,KK06},
a symmetric positive-definite matrix $P_o$ and functions $g_2$ and $g_3$ are found
to satisfy the coupled Lyapunov inequalities \eqref{coupled_lyap2} as  
$P_o=\tilde a_o\left[\begin{array}{cc}30 & -5 \\ -5&5 \end{array}\right]$
  $g_2=8\phi_{(2,3)}$, and $g_3=6\phi_{(2,3)}$,
  and with $\nu_o=5.95 \tilde a_o$, $\tilde\nu_o=18.989 \tilde a_o$, $\underline\nu_o=10 \tilde a_o$, and $\overline\nu_o=35 \tilde a_o$ with $\tilde a_o$ being any positive constant. The parameters $c_{\zeta_1}$ and $q_a$ and the initial value for $a_{\zeta_1}$ are picked as in Section~\ref{sec:example_sf}. Also, pick
  $\tilde a_o = 1$, 
  $\underline c_1 = 0.3$,
  $\underline c_2 = 0.1$, and
  $\underline c_3 = 0.0001$. 
  The initial condition for the system state vector $[x_1,x_2,x_3]^T$ is specified as $[1,-1,-1]^T$ for simulation studies. Since the initial conditions for $x_2$ and $x_3$ are not known, the initial conditions for $\hat x_2$ and $\hat x_3$ are picked simply as the values that make the initial values of the estimates for $x_2$ and $x_3$ zero, i.e., such that $\hat x_2 + rf_2(x_1)$ and $\hat x_3+r^2 f_3(x_1)$ are zero at time $t=0$. Hence, the initial condition for $[\hat x_2, \hat x_3]^T$ is $[-14.466,-21.699]^T$.  The closed-loop trajectories and the control input signal are shown in Figure~\ref{fig:sim_of}.  
  As a quick point of comparison with the ``macroscopic'' conservative approach in the controller construction in prior results (e.g., \cite{KK04f}), it can be seen that the construction of $c$ as in 
  \eqref{c_choice1} provides $c=9.791$
  while the construction of $c$ in \eqref{c_choice} provides $c=0.969$. A similar trend can be seen in each of the other design freedoms when compared between the matrix pencil based approach and the more conservative approach outlined in Section~\ref{sec:design_freedoms_of}, which is developed in more detail in \cite{KK04f}. While simulations with the design freedoms computed as in Section~\ref{sec:design_freedoms_of} are omitted here for brevity, a comparison shows that the matrix pencil based approach provides notably smoother and better transient performance due to less conservative values of design freedoms.

 \section{Conclusion}
 \label{sec:conclusion}

 A new framework was developed for dual dynamic high-gain scaling based control designs based on a matrix pencil approach that enables efficiently taking into account the detailed structure of a system. The proposed matrix pencil based design methodology significantly reduces conservatism of computed bounds in the control design thereby reducing control magnitudes and bandwidth and also reduces algebraic complexity in the computation of the design freedoms. It was shown that the proposed methodology enables efficient non-conservative computation of the design freedoms in both the state-feedback and output-feedback cases.
The efficacy of the proposed approach was demonstrated through simulation studies on an illustrative example for both the state-feedback and output-feedback scenarios.
A topic of on-going research is to determine if analogous matrix pencil based approaches can be formulated for the various other classes of systems to which dynamic high-gain scaling based control designs are applicable (e.g., feedforward systems such as \cite{KK20_matrixpencil_feedforward,KK21_matrixpencil_feedforward}). Another direction for future research is to determine whether the designs of the matrices $P_o$ and $P_c$ can be combined into an integrated methodology along with the other design freedoms instead of designing these matrices as a separate preliminary step based on coupled Lyapunov inequalities.

\bstctlcite{IEEEexample:BSTcontrol}
\bibliographystyle{IEEEtran}
\bibliography{refs}
\bstctlcite{IEEEexample:BSTcontrol}
\end{document}